\documentclass[a4paper]{article}
\usepackage{amssymb}
\usepackage{amsfonts}
\usepackage{amsmath}

\input{tcilatex}

\begin{document}

\title{The cohomology of projective unitary groups\\
{\small Dedicated to Professor Victor M. Buchstaber on his 80th birthday}}
\author{Haibao Duan \thanks{%
Supported by National Science foundation of China No. 12331003}}
\date{}
\maketitle

\begin{abstract}
The projective unitary group $PU(n)$ is the quotient of the unitary group $%
U(n)$ by its center $S^{1}=\{e^{i\theta }I_{n};\theta \in \lbrack 0,2\pi ]\}$%
, where $I_{n}$ is the identity matrix. Combining the Serre spectral
sequence of the fibration $PU(n)\rightarrow PU(n)/T$ with the Gysin sequence
of the circle bundle $U(n)\rightarrow PU(n)$, we compute the integral
cohomology ring of $PU(n)$ using explicit constructed generators, where $T$
is a maximal torus on $PU(n)$.

\begin{description}
\item[2000 Mathematical Subject Classification: ] 57T10, 55R20

\item[Key words:] Lie groups; Cohomology, Serre spectral sequence; Gysin
sequence.

\item[Email] dhb@math.ac.cn
\end{description}
\end{abstract}

\section{Introduction}

The projective unitary group $PU(n)$ is the quotient of the unitary group $%
U(n)$ by its center $S^{1}=\{e^{i\theta }I_{n};\theta \in \lbrack 0,2\pi ]\}$%
, where $I_{n}$ is the identity matrix. This is the group of holomorphic
isometries on the projective space $\mathbb{C}P^{n-1}$, essential to the
pure Yang--Mills gauge theory, and to the twisted $K$-theory \cite{A,Gr,MR}.
In this paper we compute the integral cohomology ring $H^{\ast }(PU(n))$,
using generators constructed from the Weyl invariants of the group $U(n)$.
This section gives an introduction to the main results.

In \cite{B} Borel showed that there exist elements $\xi _{2k-1}\in
H^{2k-1}(U(n))$, $1\leq k\leq n$, so that the cohomology ring $H^{\ast
}(U(n))$ is isomorphic to the exterior ring

\begin{enumerate}
\item[(1.1)] $H^{\ast }(U(n))=\Lambda (\xi _{1},\cdots ,\xi _{2n-1})$.
\end{enumerate}

\noindent In particular, as a graded group the cohomology $H^{\ast }(U(n))$
has the additive basis $\left\{ 1,\xi _{I}\right\} $, where $\xi _{I}:=%
\underset{j\in I}{\cup }${\small \ }$\xi _{2j-1}$, $I\subseteq \left\{
1,\cdots ,n\right\} $. On the other hand, regard the quotient homomorphism

\begin{enumerate}
\item[(1.2)] $c:U(n)\rightarrow PU(n)=U(n)/S^{1}$
\end{enumerate}

\noindent as an oriented circle bundle on $PU(n)$, and let $\omega \in
H^{2}(PU(n))$ be its Euler class. Then the Gysin sequence \cite[p.149]{MS}
of the oriented circle bundle $c$ gives an exact sequence, relating the
cohomologies of $PU(n)$ and $U(n)$,

\begin{enumerate}
\item[(1.3)] ${\small \cdots \rightarrow H}^{r}{\small (PU(n))}\overset{%
c^{\ast }}{\rightarrow }{\small H}^{r}{\small (U(n))}\overset{\theta }{%
\rightarrow }{\small H}^{r-1}{\small (PU(n))}\overset{\omega \cup }{%
\rightarrow }{\small H}^{r+1}{\small (PU(n))\rightarrow \cdots ,}${\small \ }
\end{enumerate}

\noindent Our aim is to solve the ring $H^{\ast }(PU(n))$ out of this
sequence. To this end we are bound to specify a set of generators of the
ring $H^{\ast }(PU(n))$, so that the induced map $c^{\ast }$ on $H^{\ast
}(PU(n))$, as well as the operator $\theta $ on $H^{\ast }(U(n))$, can be
explicitly expressed.

By \textsl{the prime factorization }of an integer\textsl{\ }$n$, we mean the
unique decomposition $n=p_{1}^{r_{1}}\cdots p_{t}^{r_{t}}$, where $%
1<p_{1}<\cdots <p_{t}$ is the set of all prime factors of $n$. It furnishes
the set $\left\{ 2,\cdots ,n\right\} $ with the partition

\begin{quote}
$\left\{ 2,\cdots ,n\right\} =Q_{0}(n)\underset{1\leq i\leq t}{\amalg }%
Q_{p_{i}}(n)$, $\quad Q_{p_{i}}(n):=\{p_{i},\cdots ,p_{i}^{r_{i}}\}$,
\end{quote}

\noindent where $Q_{0}(n)$ denotes the complement of the disjoint union $%
\underset{1\leq i\leq t}{\amalg }Q_{p_{i}}(n)$ in $\left\{ 2,\cdots
,n\right\} $. In addition, for a subset $\{a_{1},\cdots ,a_{r}\}$ of a ring $%
A$ denote by $\left\langle a_{1},\cdots ,a_{r}\right\rangle $ the ideal
generated by $a_{1},\cdots ,a_{r}$, and write $A/\left\langle a_{1},\cdots
,a_{r}\right\rangle $ for the quotient ring.

\bigskip

\noindent \textbf{Theorem A.} \textsl{There exist elements }$\rho _{2k-1}\in
H^{2k-1}(PU(n))$\textsl{,} $2\leq k\leq n$\textsl{,\ such that}

\begin{enumerate}
\item[(1.4)] $c^{\ast }(\rho _{2k-1})=\left\{ 
\begin{array}{c}
\xi _{2k-1}\text{ \textsl{if} }k\in Q_{0}(n)\text{\textsl{,}} \\ 
p_{i}\cdot \xi _{2k-1}\text{ \textsl{if} }k\in Q_{p_{i}}(n)\text{\textsl{.}}%
\end{array}%
\right. $
\end{enumerate}

\textsl{Moreover,} \textsl{the ring }$H^{\ast }(PU(n))$\textsl{\ has the
presentation}

\begin{enumerate}
\item[(1.5)] $H^{\ast }(PU(n))=\frac{\mathbb{Z}[\omega ]\otimes \Lambda
(\rho _{3},\cdots ,\rho _{2n-1})}{\langle \omega \otimes \theta (\xi
_{I}),I\subseteq \{1,\cdots ,n\}\rangle }$\textsl{,}
\end{enumerate}

\noindent \textsl{where }$\theta :H^{\ast }(U(n))\rightarrow H^{\ast
}(PU(n)) $ \textsl{is the homomorphism} \textsl{in (1.3).}

\bigskip

In the view of formula (1.5), the determination of the ring $H^{\ast
}(PU(n)) $ amounts to find the formula that expresses the terms $\theta (\xi
_{I})$ as elements in the ring $\mathbb{Z}[\omega ]\otimes \Lambda (\rho
_{3},\cdots ,\rho _{2n-1})$. Our second result serves this purpose.

\bigskip

\noindent \textbf{Theorem B. }\textsl{The elements }$\theta (\xi _{I})\in 
\mathbb{Z}[\omega ]\otimes \Lambda (\rho _{3},\cdots ,\rho _{2n-1})$ \textsl{%
in (1.5)} \textsl{are given by the following recurrence relations:}

\begin{enumerate}
\item[a)] $\theta (\xi _{2k-1})=\binom{n}{k}\omega ^{k-1}$\textsl{, }$k\in
\left\{ 1,\cdots ,n\right\} $\textsl{;}

\item[b)] $\theta (\xi _{I})=-\theta (\xi _{I^{e}})\cup \rho _{2i_{k}-1}$ 
\textsl{if }$i_{k}\in Q_{0}(n)$\textsl{;}

\item[c)] $\theta (\xi _{I})=-\frac{1}{p_{i}}(\theta (\xi _{I^{e}})\cup \rho
_{2i_{k}-1}+\omega ^{i_{k}-\frac{i_{k}}{p_{i}}}\cup \theta (\xi _{I^{e}}\cup
\xi _{2\frac{i_{k}}{p_{i}}-1})$ \textsl{if }$i_{k}\in Q_{p_{i}}(n)$\textsl{,}
\end{enumerate}

\noindent \textsl{where }$I=\left\{ i_{1},\cdots ,i_{k}\right\} \subset
\left\{ 1,\cdots ,n\right\} $\textsl{\ is a multi-index with} \textsl{length}
$k\geq 2$\textsl{, and where }$I^{e}:=${\small \ }$\{i_{1},\cdots ,i_{k-1},%
\widehat{i_{k}}\}$.

\bigskip

Granted with Theorems A and B our main result, stated in Theorem C of
Section \S 6, gives a concise and useful presentation of the ring $H^{\ast
}(PU(n))$.

The study of the cohomology of the group $PU(n)$ has a long history.
Notably, Borel \cite{B1} computed the $\func{mod}p$ cohomology ring $H^{\ast
}(PU(n);\mathbb{F}_{p})$ for $p$ a prime. In \cite{BB} Baum and Browder
recovered Borel's calculation, where an approach to the integral computation
is sketched, but not completed. In \cite{R} Ruiz gave a presentation of the
ring $H^{\ast }(PU(n))$ which, as pointed out by Astey, Gitler, Micha and
Pastor \cite{AGMP}, contains unrecoverable errors. In the course of studying
topological period-index problem, Antieau and Williams \cite{AW1} computed
the ring $H^{\ast }(PU(n))$ up to $n=5$. Another relevant work is due to
Petrie \cite{P1,P2}: he proved that, if $n$ is a prime power $p^{r}$, then $%
H^{\ast }(PU(n))$ is additively isomorphic to the topological K-theory $%
K(PU(n))$. In a recent paper \cite{C} Can pointed out that our result gives
also the cohomology ring of the determinantal variety of the monoid of $%
n\times n$ matrices space.

The main tool in our construction and computation is the Serre spectral
sequence $\left\{ E_{r}^{\ast ,\ast }(G),d_{r}\right\} $ of the torus
fibration

\begin{quote}
$T\rightarrow G\rightarrow G/T$,
\end{quote}

\noindent where $G=U(n)$ or $PU(n)$, $T$ is a fixed maximal torus on $G$. In
Section 3 we make the second page $\left\{ E_{2}^{\ast ,\ast
}(G),d_{2}\right\} $ precise for both $G=U(n)$ and $PU(n)$. In Section 4 we
establish an exact sequence relating the third page $E_{3}^{\ast ,\ast
}(PU(n))$ with $E_{3}^{\ast ,\ast }(U(n))$, and construct the preliminary
form $\rho _{2k-1}^{\prime }\in E_{3}^{2k-2,1}(U(n))$ for the cohomology
class $\rho _{2k-1}\in H^{2k-1}(PU(n))$ promised in Theorem A. After solving
the extension problem from $E_{3}^{\ast ,\ast }(PU(n))$ to the cohomology $%
H^{\ast }(PU(n))$, Theorems A and B are proved in Section 5. Finally,
Section \S 6 is devoted to show the main result Theorem C, as well as its
applications

\bigskip

\noindent \textbf{Acknowledgement.} The author would like to thank his
referee for improvements on the early version of the present paper.

\section{The binomial coefficients $\binom{n}{r}$}

Our construction of the ring $H^{\ast }(PU(n))$ relies on two arithmetic
results about the binomial coefficients $\binom{n}{r}=\frac{n!}{r!(n-r)!}$.
For each $1\leq k\leq n$ consider the integer

\begin{quote}
$b_{n,k}:=g.c.d.\{\binom{n}{1},\cdots ,\binom{n}{k}\}$.
\end{quote}

\noindent Since $b_{n,k}$ divides $b_{n,k-1}$ the ratio $\frac{b_{n,k-1}}{%
b_{n,k}}$ is always an integer. The following result, proved in \cite{DL},
expresses the integer $b_{n,k}$ as a function in $k$ and $n$.

\bigskip

\noindent \textbf{Lemma 2.1. }\textsl{If the integer }$n$\textsl{\ has the
prime factorization }$p_{1}^{r_{1}}\cdots p_{t}^{r_{t}}$\textsl{, then}

\begin{quote}
\textsl{i) }$\frac{b_{n,k-1}}{b_{n,k}}=\left\{ 
\begin{array}{c}
1\text{ \textsl{if} }k\in Q_{0}(n)\text{\textsl{,}} \\ 
p_{i}\text{ \textsl{if} }k\in Q_{p_{i}}(n)\text{\textsl{.}}%
\end{array}%
\right. $

\textsl{ii) }$b_{n,k}=p_{1}^{\varepsilon _{1}}\cdots p_{t}^{\varepsilon
_{t}} $\textsl{, }
\end{quote}

\noindent \textsl{where }$\varepsilon _{i}$\textsl{\ is the cardinality of
the set of intersection} $\left\{ k,\cdots ,n\right\} \cap Q_{p_{i}}(n)$%
\textsl{.}\hfill $\square $

\bigskip

Let $e_{k}(x_{1},\cdots ,x_{n})$\ be the $k^{th}$\ elementary symmetric
function in $x_{1},\cdots ,x_{n}$, and let $g_{k}(x_{1},\cdots
,x_{n})=x_{1}^{k}+\cdots +x_{n}^{k}$ be the $k^{th}$ power sum. In the
Newton formula

\begin{quote}
$k\cdot e_{k}=\underset{1\leq i\leq k}{\Sigma }(-1)^{i-1}g_{i}e_{k-i}$,
\end{quote}

\noindent taking $x_{1}=\cdots =x_{n}=1$, one gets from

\begin{quote}
$g_{i}(1,\cdots ,1)=n$ and $e_{k-i}(1,\cdots ,1)=\binom{n}{k-i}$
\end{quote}

\noindent the expansion

\begin{enumerate}
\item[(2.1)] $\binom{n}{k}=\frac{n}{k}\cdot \underset{1\leq i\leq k}{\Sigma }%
(-1)^{i-1}\binom{n}{k-i}$.
\end{enumerate}

\bigskip

\noindent \textbf{Lemma 2.2. }\textsl{Suppose that }$p$\textsl{\ is a prime, 
}$n=p^{r}\cdot n^{\prime }$\textsl{\ with }$(p,n^{\prime })=1$\textsl{.\ Then%
}

\begin{enumerate}
\item[(2.2)] $\binom{n}{p^{s}}=\frac{n}{p^{s}}\cdot \underset{1\leq i\leq
p^{s}-p^{s-1}}{\Sigma }(-1)^{i-1}\binom{n}{p^{s}-i}+\frac{1}{p}\binom{n}{%
p^{s-1}}$\textsl{,} $1\leq s\leq r$\textsl{.}
\end{enumerate}

\noindent \textbf{Proof.} For $k=p^{s}$ separate the sum (2.1) at $%
i=p^{s}-p^{s-1}$ to get

\begin{quote}
$\binom{n}{p^{s}}=\frac{n}{p^{s}}\cdot \underset{1\leq i\leq p^{s}-p^{s-1}}{%
\Sigma }(-1)^{i-1}\binom{n}{p^{s}-i}+\frac{n}{p^{s}}\cdot \underset{%
p^{s}-p^{s-1}+1\leq i\leq p^{s}}{\Sigma }(-1)^{i-1}\binom{n}{p^{s}-i}$.
\end{quote}

\noindent Introducing $i=t+p^{s}-p^{s-1}$ the second sum on the right hand
side becomes

\begin{quote}
$(-1)^{p^{s}-p^{s-1}}\frac{n}{p^{s}}\cdot \underset{1\leq t\leq p^{s-1}}{%
\Sigma }(-1)^{t-1}\binom{n}{p^{s-1}-t}=\frac{(-1)^{p^{s}-p^{s-1}}}{p}\cdot 
\binom{n}{p^{s-1}}$
\end{quote}

\noindent by (2.1). Since $p$ is a prime, we get (2.2) from $%
(-1)^{p^{s}-p^{s-1}}=1$.\hfill $\square $

\section{Spectral sequence of the fibration $\protect\pi :$ $G\rightarrow
G/T $}

For a compact and connected Lie group $G$ with a maximal torus $T$ the
inclusion $i:T\rightarrow G$ induces the sequence of fibrations

\begin{enumerate}
\item[(3.1)] $T\rightarrow G\overset{\pi }{\rightarrow }G/T\overset{f}{%
\rightarrow }B_{T}\overset{Bi}{\rightarrow }B_{G}$,
\end{enumerate}

\noindent where $B_{H}$ denotes the classifying space of a Lie group $H$.
The quotient space $G/T$ is a smooth projective variety, called \textsl{the
complete flag manifold} of $G$ \cite{DZ}. Let $\left\{ E_{r}^{\ast ,\ast
}(G),d_{r}\right\} $ be the Serre spectral sequence of the fibration $\pi $,
and let $\mathcal{F}^{p}$ be the filtration on $H^{\ast }(G)$ defined by $%
\pi $. Then (i.e. \cite[P.146]{Mc})

\begin{center}
$0=\mathcal{F}^{r+1}(H^{r}(G))\subseteq\mathcal{F}^{r}(H^{r}(G))\subseteq
\cdots\subseteq\mathcal{F}^{0}(H^{r}(G))=H^{r}(G)$
\end{center}

\noindent with

\begin{enumerate}
\item[(3.2)] $E_{\infty }^{p,q}(G)=\mathcal{F}^{p}(H^{p+q}(G))/\mathcal{F}%
^{p+1}(H^{p+q}(G))$.
\end{enumerate}

In this section we recall basic properties of the second page $\{E_{2}^{\ast
,\ast }(G),d_{2}\}$ of the Serre spectral sequence of the map $\pi $. For
the cases $G=PU(n)$ and $U(n)$ of our interests, we bring a connection from $%
E_{2}^{\ast ,\ast }(PU(n))$ to $E_{2}^{\ast ,\ast }(U(n))$, and clarify two
canonical maps

\begin{quote}
$\pi ^{\ast }:E_{3}^{r,0}(G)\rightarrow \mathcal{F}^{r}(H^{r}(G))\subseteq
H^{r}(G)$, and

$\kappa :E_{3}^{2s,1}(G)\rightarrow \mathcal{F}^{2s}(H^{2s+1}(G))\subset
H^{\ast }(G)$,
\end{quote}

\noindent useful for solving the extension problem from $E_{3}^{\ast ,\ast
}(G)$ to $H^{\ast }(G)$.

\subsection{The complex $\{E_{2}^{\ast ,\ast }(G),d_{2}\}$ for $G=U(n)$ or $%
PU(n)$}

Recall from \cite{D1} that the\textsl{\ Borel} \textsl{transgression}\textbf{%
\ }in the fibration $\pi $ is the composition

\begin{quote}
$\tau =$ $(\pi ^{\ast })^{-1}\circ \delta :H^{1}(T)\overset{\delta }{%
\rightarrow }H^{2}(G,T)\overset{(\pi ^{\ast })^{-1}}{\rightarrow }H^{2}(G/T)$%
,
\end{quote}

\noindent where $\delta $ is the connecting homomorphism in the
cohomological exact sequence of the pair $(G,T)$, and where the map $\pi
^{\ast }$ from $H^{2}(G/T)$ to $H^{2}(G,T)$ is always an isomorphism.

\bigskip

\noindent \textbf{Lemma 3.1.} \textsl{The complex} $\left\{ E_{2}^{\ast
,\ast }(G),d_{2}\right\} $ \textsl{is}

\begin{enumerate}
\item[(3.3)] $E_{2}^{\ast ,\ast }(G)=H^{\ast }(G/T)\otimes H^{\ast }(T)$
\end{enumerate}

\noindent \textsl{on which }$d_{2}$ \textsl{is the antiderivation} \textsl{%
of degree} $1$ \textsl{given by the transgression }$\tau $ \textsl{as}

\begin{quote}
$d_{2}(x\otimes 1)=0$\textsl{,} $x\in H^{\ast }(G/T)$;

$d_{2}(1\otimes t)=\tau (t)\otimes 1$\textsl{, }$t\in H^{1}(T)$.
\end{quote}

\noindent \textsl{In addition, if }$\dim G/T=m$\textsl{,} $\dim T=n$\textsl{%
, then}

\begin{enumerate}
\item[(3.4)] $E_{r}^{m,n}(G)=H^{n+m}(G)=\mathbb{Z}$ \textsl{for all} $r\geq
2 $\textsl{;}

\item[(3.5)] $\dim E_{3}^{\ast ,\ast }(G)\otimes \mathbb{R}=\dim H^{\ast }(G;%
\mathbb{R})\mathbb{=}2^{n}$\textsl{;}

\item[(3.6)] $\rho ^{2}=0$ \textsl{for any} $\rho \in E_{3}^{2k,1}(G)$%
\textsl{;}

\item[(3.7)] $E_{3}^{\ast ,0}(G)=H^{\ast }(G/T)/\left\langle \func{Im}\tau
\right\rangle $\textsl{,}
\end{enumerate}

\noindent \textsl{where }$\left\langle \func{Im}\tau \right\rangle $\textsl{%
\ is the ideal of ring }$H^{\ast }(G/T)$ \textsl{generated by the subgroup} $%
\func{Im}\tau $\textsl{.}

\bigskip

\noindent \textbf{Proof.} Since the base manifold $G/T$ of $\pi $ is simply
connected \cite{BS}, the formula (3.3) of $E_{2}^{\ast ,\ast }(G)$ follows
from Leray-Serre theorem \cite[p.135]{Mc}. With respect to the presentation
(3.3), the $d_{2}$ action on $E_{2}^{\ast ,\ast }(G)$ is also known \cite[%
Formula (1.2)]{D1}.

Since both $G/T$ and $T$ are smooth and orientable manifolds with dimensions 
$m$ and $n$, respectively, we get by (3.3) that

\begin{quote}
$E_{2}^{m,n}(G)=\mathbb{Z}$, $E_{2}^{s,t}(G)=0$ for either $s>m$ or $t>n$.
\end{quote}

\noindent These imply that any differential $d_{r}$ that acts or lands on $%
E_{r}^{m,n}(G)$ must be trivial, showing the isomorphisms in (3.4):

\begin{quote}
$E_{2}^{m,n}(G)=E_{3}^{m,n}(G)=\cdots =E_{\infty }^{m,n}(G)=H^{n+m}(G)$.
\end{quote}

With $E_{3}^{\ast ,\ast }(G)\otimes \mathbb{R=}H^{\ast }(G;\mathbb{R})$ by
Leray \cite{L} we obtain (3.5) by the classical result of Hopf \cite{Ho}
that $H^{\ast }(G;\mathbb{R})$ is an exterior ring in $n$ generators.

For a $d_{2}$-cocycle $\xi \in E_{2}^{2k,1}(G)$ write $[\xi ]\in E_{3}^{\ast
,1}(G)$ for its cohomology class. One gets $[\xi ]^{2}=0$ from $\xi ^{2}=0$
by the multiplicative rule $t_{1}t_{2}=-t_{2}t_{1}$ for $t_{1},t_{2}\in
H^{1}(T)$, showing property (3.6).

Finally, since $E_{3}^{\ast ,0}(G)$ is the cokernel of the differential

\begin{quote}
$d_{2}:E_{2}^{\ast ,1}(G)=H^{\ast }(G/T)\otimes H^{1}(T)\rightarrow
E_{2}^{\ast ,0}(G)=H^{\ast }(G/T)$,
\end{quote}

\noindent we get (3.7) from $\func{Im}d_{2}=\left\langle \func{Im}\tau
\right\rangle \subset H^{\ast }(G/T)$ by above formulae of $d_{2}$.\hfill $%
\square $

\bigskip

\noindent \textbf{Example 3.2.} For $G=U(n)$ the complexes $\left\{
E_{2}^{\ast ,\ast }(G),d_{2}\right\} $ has been determined by Borel.
Precisely, take the diagonal subgroup $diag\{e^{i\theta _{1}},\cdots
,e^{i\theta _{n}}\}$ as the fixed maximal torus $T$ on $U(n)$, and let $%
y_{i}\in H^{1}(T)$ be the Kronecker dual of the oriented circle subgroup $%
\theta _{k}=0$, $k\neq i$, on $T$. Set

\begin{quote}
$x_{i}:=\tau (y_{i})\in H^{2}(U(n)/T)$, $1\leq i\leq n$,
\end{quote}

\noindent and let $e_{r}\in \mathbb{Z}[x_{1},\cdots ,x_{n}]$ be the $r^{th}$%
\ elementary symmetric polynomial in $x_{1},\cdots ,x_{n}$.

Borel \cite{B} has shown that $H^{\ast }(T)=\Lambda (y_{1},\cdots ,y_{n})$,
and that, with respect to the canonical presentations

\begin{quote}
$H^{\ast }(B_{T})=\mathbb{Z}[x_{1},\cdots ,x_{n}]$, $H^{\ast }(U(n)/T)=\frac{%
\mathbb{Z}[x_{1},\cdots ,x_{n}]}{\left\langle e_{1},\cdots
,e_{n}\right\rangle }$,
\end{quote}

\noindent the induced map $f^{\ast }$ (see (3.1)) is the obvious quotient map

\begin{quote}
$\quad H^{\ast }(B_{T})=\mathbb{Z}[x_{1},\cdots ,x_{n}]\rightarrow H^{\ast
}(U(n)/T)=\frac{\mathbb{Z}[x_{1},\cdots ,x_{n}]}{\left\langle e_{1},\cdots
,e_{n}\right\rangle }$.
\end{quote}

\noindent According to Lemma 3.1

\begin{enumerate}
\item[(3.8)] $E_{2}^{\ast ,\ast }(U(n))=\frac{\mathbb{Z}[x_{1},\cdots ,x_{n}]%
}{\left\langle e_{1},\cdots ,e_{n}\right\rangle }\otimes \Lambda
(y_{1},\cdots ,y_{n})$ with

$d_{2}(x_{i}\otimes 1)=0$, $d_{2}(1\otimes y_{i})=x_{i}\otimes 1$.\hfill $%
\square $
\end{enumerate}

\bigskip

Turning to the group $PU(n)$ we take $T^{\prime }:=c(T)$ as the preferable
maximal torus, and denote by $\pi $ and $\pi ^{\prime }$ for the
corresponding torus fibrations over $U(n)/T$ and $PU(n)/T^{\prime }$,
respectively. Then the circle bundle $c$ on $PU(n)$ (i.e. (1.2)) can be
regarded as a bundle map from $\pi $ to $\pi ^{\prime }$:

\begin{center}
$%
\begin{array}{ccccc}
S^{1} & \overset{i}{\hookrightarrow } & T & \overset{c^{\prime }}{%
\rightarrow } & T^{\prime } \\ 
\parallel &  & \cap \quad &  & \cap \quad \\ 
S^{1} & \hookrightarrow & U(n) & \overset{c}{\rightarrow } & PU(n) \\ 
&  & \pi \downarrow \quad &  & \pi ^{\prime }\downarrow \quad \\ 
&  & U(n)/T & = & PU(n)/T^{\prime }%
\end{array}%
$,
\end{center}

\noindent where $c^{\prime }$ denotes the restriction of $c$ on the torus $T$%
, and where the base manifolds $U(n)/T$ and $PU(n)/T^{\prime }$ are
canonically isomorphic \cite{DZ}. In order to get a presentation of the ring 
$E_{2}^{\ast ,\ast }(PU(n))$ so that the induced map $c^{\ast }:E_{2}^{\ast
,\ast }(PU(n))\rightarrow E_{2}^{\ast ,\ast }(U(n))$ is transparent, we
change the basis $\left\{ y_{1},\cdots ,y_{n}\right\} $ of $H^{1}(T)$ by
setting

\begin{quote}
$t_{0}:=y_{1}$, $t_{i}:=y_{i+1}-y_{i}$ , $1\leq i\leq n-1$.
\end{quote}

\noindent \textbf{Lemma 3.3.} \textsl{With the notation above} \textsl{one
has}

\begin{quote}
$E_{2}^{\ast ,\ast }(U(n))=H^{\ast }(U(n)/T)\otimes \Lambda
(t_{0},t_{1},\cdots ,t_{n-1})$\textsl{;}

$E_{2}^{\ast ,\ast }(PU(n))=H^{\ast }(U(n)/T)\otimes \Lambda (t_{1},\cdots
,t_{n-1})$.
\end{quote}

\noindent \textsl{Moreover}

\begin{quote}
\textsl{i)} \textsl{The map} $c^{\ast }:$\textsl{\ }$E_{2}^{\ast ,\ast
}(PU(n))\rightarrow E_{2}^{\ast ,\ast }(U(n))$ \textsl{is the inclusion}

$\qquad c^{\ast }(z\otimes t_{i})=z\otimes t_{i}$\textsl{, }$z\in H^{\ast
}(U(n)/T)$\textsl{, }$1\leq i\leq n-1;$

\textsl{ii) The differential }$d_{2}$ \textsl{on} $E_{2}^{\ast ,\ast }(U(n))$%
\textsl{\ (resp. }$d_{2}^{\prime }$ \textsl{on} $E_{2}^{\ast ,\ast }(PU(n))$%
\textsl{)} \textsl{is}

$\qquad d_{2}(1\otimes t_{0})=x_{1}\otimes 1$\textsl{,}

$\qquad d_{2}(1\otimes t_{i})$ $=d_{2}^{\prime }(1\otimes
t_{i})=(x_{i+1}-x_{i})\otimes 1$\textsl{, for }$1\leq i\leq n-1$\textsl{.}
\end{quote}

\textsl{In particular, letting }$\overline{\omega }\in E_{3}^{2,0}(PU(n))$%
\textsl{\ be the cohomology class }$[x_{1}\otimes 1]$ \textsl{of the }$%
d_{2}^{\prime }$\textsl{-cocycle} $x_{1}\otimes 1\in E_{2}^{2,0}(PU(n))$%
\textsl{, then}

\begin{enumerate}
\item[(3.9)] $E_{3}^{\ast ,0}(PU(n))=\frac{\mathbb{Z}[\overline{\omega }]}{%
\left\langle b_{n,k}\overline{\omega }^{k},1\leq k\leq n\right\rangle }$%
\textsl{,\ where }$b_{n,k}=g.c.d.\{\binom{n}{1},\cdots ,\binom{n}{k}\}$%
\textsl{.}
\end{enumerate}

\noindent \textbf{Proof.} The maximal torus $T^{\prime }$ of $PU(n)$ has the
factorization $S^{1}\times \cdots \times S^{1}$ ($n-1$ factors) so that the
map $c^{\prime }:T\rightarrow T^{\prime }$ is given by

\begin{quote}
$c^{\prime }(diag\{e^{i\theta _{1}},\cdots ,e^{i\theta _{n}}\})=(e^{i(\theta
_{2}-\theta _{1})},\cdots ,e^{i(\theta _{n}-\theta _{n-1})})$.
\end{quote}

\noindent In particular, $H^{\ast }(T)=\Lambda (t_{0},t_{1},\cdots ,t_{n-1})$%
, while the map $c^{\prime \ast }$ identifies $H^{\ast }(T^{\prime })$ with
the subring $\Lambda (t_{1},\cdots ,t_{n-1})$ of $H^{\ast }(T)$. This
justifies the formulae of $E_{2}^{\ast ,\ast }(U(n))$ and $E_{2}^{\ast ,\ast
}(PU(n))$ stated in the lemma.

Since $c^{\ast }$ is a map of $d_{2}$-complexes, the formulae i) and ii) are
clear by (3.8). For the expression (3.9) of the subring $E_{3}^{\ast
,0}(PU(n))\subset E_{3}^{\ast ,\ast }(PU(n))$ we note that $\func{Im}\tau
^{\prime }$ has the basis $\{x_{2}-x_{1},\cdots ,x_{n}-x_{n-1}\}$ by the
formula of $d_{2}^{\prime }$. In particular, with

\begin{quote}
$e_{k}(x_{1},\cdots ,x_{n})\equiv \binom{n}{k}x_{1}^{k}\func{mod}\func{Im}%
\tau ^{\prime }$
\end{quote}

\noindent we get from (3.7) that

\begin{quote}
$E_{3}^{\ast ,0}(PU(n))=\frac{\mathbb{Z}[x_{1},\cdots ,x_{n}]}{\left\langle
e_{1},\cdots ,e_{n}\right\rangle }/\left\langle \func{Im}\tau ^{\prime
}\right\rangle =\frac{\mathbb{Z}[x_{1}]}{\left\langle \binom{n}{k}%
x_{1}^{k},1\leq k\leq n\right\rangle }$.
\end{quote}

\noindent This is equivalent to (3.9), because the relations $\binom{n}{k}%
x_{1}^{k}=0$ on $E_{3}^{\ast ,0}(PU(n))$, $1\leq k\leq n$, imply that the
order of the power $\overline{\omega }^{r}$ is precisely $b_{n,r}$.\hfill $%
\square $

\subsection{The map $\protect\pi ^{\ast }:E_{3}^{\ast ,0}(G)\rightarrow
H^{\ast }(G)$}

Since the differential $d_{r}$ on $E_{r}^{p,q}(G)$ is of bi-degree $(r,-r+1)$%
, we must have $d_{r}(E_{r}^{\ast ,0}(G))=0$ for all $r\geq 2$. This gives
rise to a sequence of quotient maps

\begin{center}
$H^{r}(G/T)=E_{2}^{r,0}\twoheadrightarrow E_{3}^{r,0}\twoheadrightarrow
\cdots\twoheadrightarrow E_{\infty}^{r,0}=\mathcal{F}^{r}(H^{r}(G))\subseteq
H^{r}(G)$
\end{center}

\noindent whose composition agrees with the induced map $\pi ^{\ast
}:H^{\ast }(G/T)\rightarrow H^{\ast }(G)$ (\cite[P.147]{Mc}). Therefore, we
may reserve $\pi ^{\ast }$ for the epimorphism:

\begin{enumerate}
\item[(3.10)] $\pi ^{\ast }:E_{3}^{\ast ,0}(G)\rightarrow \func{Im}\pi
^{\ast }\subset H^{\ast }(G)$.
\end{enumerate}

\noindent \textbf{Lemma 3.4.} \textsl{Let }$\omega \in H^{2}(PU(n))$\textsl{%
\ be the Euler class of the circle bundle }$c$\textsl{.} \textsl{Then the
map }$\pi ^{\prime }$ \textsl{satisfies, in terms of (3.9) and (3.10), the
following relation}

\begin{enumerate}
\item[(3.11)] $\pi ^{\prime \ast }(\overline{\omega })=$ $\omega $.
\end{enumerate}

\noindent \textbf{Proof. }Let\textbf{\ }$i:SU(n)\rightarrow U(n)$ be the
inclusion of the special unitary group $SU(n)$. Then the composition $%
\widetilde{c}=c\circ i:$ $SU(n)\rightarrow PU(n)$ is the universal covering
of $PU(n)$ with $\deg \widetilde{c}=n$. From the homotopy exact sequence of $%
\widetilde{c}$ one finds that

\begin{quote}
$\pi _{1}(PU(n))=\mathbb{Z}_{n}$ and $\pi _{2}(PU(n))=0$,
\end{quote}

\noindent implying that $H^{1}(PU(n))=0$ and $H^{2}(PU(n))=\mathbb{Z}_{n}$,
in which $\omega $ is a generator of $H^{2}(PU(n))=\mathbb{Z}_{n}$ by (1.3)
(e.g. \cite{DL}).

Consider the part of the cohomology exact sequence of the pair $%
(PU(n),T^{\prime })$

\begin{center}
$%
\begin{array}{cccccc}
H^{1}(PU(n)) & \rightarrow & H^{1}(T^{\prime }) & \overset{\delta }{%
\rightarrow }H^{2}(PU(n),T^{\prime }) & \rightarrow & H^{2}(PU(n))=\mathbb{Z}%
_{n} \\ 
\parallel &  & \tau ^{^{\prime }}\searrow & (\pi ^{\prime \ast
})^{-1}\downarrow \cong & \nearrow \pi ^{\prime \ast } &  \\ 
0 &  &  & H^{2}(PU(n)/T^{\prime }) &  & 
\end{array}%
$,
\end{center}

\noindent where $PU(n)/T^{\prime }=U(n)/T$, and where the composition $(\pi
^{\prime \ast })^{-1}\circ \delta $ is the transgression $\tau ^{\prime }$
in $\pi ^{\prime }$ which, by Lemma 3.1 and ii) of Lemma 3.3, is given by

\begin{quote}
$\tau ^{\prime }(t_{i})=x_{i+1}-x_{i}$, $1\leq i\leq n-1$.
\end{quote}

\noindent It follows that the map $\pi ^{\prime \ast }$ induces a
monomorphism

\begin{quote}
$H^{2}(PU(n)/T^{\prime })/\func{Im}\tau ^{\prime }\rightarrow H^{2}(PU(n))=%
\mathbb{Z}_{n}$.
\end{quote}

\noindent We obtain (3.11) from $E_{3}^{2,0}(PU(n))=H^{2}(PU(n)/T^{\prime })/%
\func{Im}\tau ^{\prime }$, which is a cyclic group of order $n$ with
generator $\overline{\omega }$ by (3.9). \hfill \hfill $\square $

\subsection{The map $\protect\kappa :E_{3}^{\ast ,1}(G)\rightarrow H^{\ast
}(G)$}

Since the group $E_{r}^{p,q}(G)$ is a subquotient of $E_{r-1}^{p,q}(G)$, the
fact $H^{2s+1}(G/T)=0$ by Bott-Samelson \cite{BS} implies that

\begin{quote}
$E_{r}^{2s+1,q}(G)=0$ for all $s,q\geq 0$ and $r\geq 2$.
\end{quote}

\noindent In particular, taking $(r,q)=(\infty ,0)$ one gets by (3.2) that

\begin{quote}
$\mathcal{F}^{2s+1}(H^{2s+1}(G))=\mathcal{F}^{2s+2}(H^{2s+1}(G))=0$.
\end{quote}

\noindent It implies, again by (3.2), that

\begin{quote}
$E_{\infty }^{2s,1}(G)=\mathcal{F}^{2s}(H^{2s+1}(G))\subset H^{2s+1}(G)$.
\end{quote}

\noindent Moreover, Since the differential $d_{r}$ on $E_{r}^{p,q}(G)$ is of
bi-degree $(r,r-1)$, we get $d_{r}(E_{r}^{\ast ,1})=0$ for all $r\geq 3$.
Summarizing, we obtain the sequence of epimorphisms

\begin{enumerate}
\item[(3.12)] $\kappa :E_{3}^{\ast ,1}(G)\twoheadrightarrow E_{4}^{\ast
,1}(G)\twoheadrightarrow \cdots \twoheadrightarrow E_{\infty }^{\ast ,1}(G)=%
\mathcal{F}^{2s}(H^{2s+1}(G))\subset H^{\ast }(G)$
\end{enumerate}

\noindent that interprets elements of $E_{3}^{2k,1}(G)$ directly as
cohomology classes of $G$.

The map $\kappa $ is essential for us to solve the extension problem from $%
E_{3}^{\ast ,\ast }(G)$ to the cohomology $H^{\ast }(G)$, where $G=U(n)$ and 
$PU(n)$. As an initial step we examine the case $G=U(n)$ for which

\begin{quote}
$E_{2}^{\ast ,\ast }(U(n))=\frac{\mathbb{Z}[x_{1},\cdots ,x_{n}]}{%
\left\langle e_{1},\cdots ,e_{n}\right\rangle }\otimes \Lambda (t_{0},\cdots
,t_{n-1})$ (by Lemma 3.3).
\end{quote}

\noindent In term of the symmetric polynomials $e_{k}\in \mathbb{Z}%
[x_{1},\cdots ,x_{n}]$ therein we will now define $d_{2}$-cocycles $\widehat{%
e}_{k}\in E_{2}^{2k-2,1}(U(n))$ whose cohomology classes generate the rings $%
E_{3}^{\ast ,\ast }(U(n))$.

In general, by \textsl{the Taylor expansion }of a homogeneous polynomial $%
h\in H^{\ast }(B_{T})=\mathbb{Z}[x_{1},\cdots ,x_{n}]$ we mean the unique
expression

\begin{enumerate}
\item[(3.13)] $h=h^{(1)}\cdot x_{1}+h^{(2)}\cdot (x_{2}-x_{1})+\cdots
+h^{(n)}\cdot (x_{n}-x_{n-1})$\textsl{,}
\end{enumerate}

\noindent such that\textsl{\ }$h^{(r)}\in \mathbb{Z}[x_{1},\cdots ,x_{r}]$.%
\textsl{\ }As examples, we have

\begin{enumerate}
\item[(3.14)] $h^{(1)}=\frac{h(x_{1},\cdots ,x_{1})}{x_{1}}\in \mathbb{Z}%
[x_{1}]$;

$h^{(2)}=\frac{h(x_{1},x_{2},\cdots ,x_{2})-h(x_{1},\cdots ,x_{1})}{%
x_{2}-x_{1}}\in \mathbb{Z}[x_{1},x_{2}]$, $\cdots $, etc.
\end{enumerate}

\noindent In view of the expansion (3.13) define the linear map

\begin{quote}
$\widehat{}:\mathbb{Z}[x_{1},\cdots ,x_{n}]\rightarrow E_{2}^{\ast ,1}(U(n))$
\end{quote}

\noindent by letting

\begin{enumerate}
\item[(3.15)] $\widehat{h}:=h^{(1)}\otimes t_{0}+h^{(2)}\otimes t_{1}+\cdots
+h^{(n)}\otimes t_{n-1}\in E_{2}^{\ast ,1}(U(n))$.
\end{enumerate}

\noindent By ii) of Lemma 3.2, the map\ $\widehat{}$ is a lift of the
quotient map $f^{\ast }$ relative to $d_{2}$:

\begin{enumerate}
\item[(3.16)] $d_{2}(\widehat{h})=f^{\ast }(h):$%
\begin{tabular}{lll}
&  & $E_{2}^{\ast ,1}(U(n))$ \\ 
& $\widehat{}\nearrow $ & $d_{2}\downarrow $ \\ 
$\mathbb{Z}[x_{1},\cdots ,x_{n}]$ & $\overset{f^{\ast }}{\longrightarrow }$
& $E_{2}^{\ast ,0}(U(n))=\frac{\mathbb{Z}[x_{1},\cdots ,x_{n}]}{\left\langle
e_{1},\cdots ,e_{n}\right\rangle }$.%
\end{tabular}
\end{enumerate}

\bigskip

\noindent \textbf{Definition 3.5.} For a $d_{2}$-cocycle $z\in E_{2}^{\ast
,\ast }(G)$ write $\left[ z\right] \in E_{3}^{\ast ,\ast }(G)$ to denote its
cohomology class.

For the symmetric polynomial $e_{k}\in \mathbb{Z}[x_{1},\cdots ,x_{n}]$ the
lift $\widehat{e}_{k}\in E_{2}^{2k-2,1}(U(n))$ is a $d_{2}$-cocycle by
(3.16):

\begin{quote}
$d_{2}(\widehat{e_{k}})=f^{\ast }(e_{k})=0$.
\end{quote}

\noindent Define, successively, the cohomology classes $\xi _{2k-1}^{\prime
} $ and $\xi _{2k-1}$ by

\begin{enumerate}
\item[(3.17)] $\xi _{2k-1}^{\prime }:=\left[ \widehat{e}_{k}\right] \in
E_{3}^{2k-2,1}(U(n))$, and

$\xi _{2k-1}:=\kappa (\xi _{2k-1}^{\prime })\in H^{2k-1}(U(n))$, $1\leq
k\leq n$.
\end{enumerate}

\noindent To emphasize that the class $\xi _{2k-1}$ is defined in sequel to $%
\xi _{2k-1}^{\prime }$, we call $\xi _{2k-1}^{\prime }$ \textsl{the} \textsl{%
primary }$1$\textsl{-form} of $\xi _{2k-1}$.\hfill $\square $

\bigskip

It can be shown that (e.g. Pittie \cite[4.2.Proposition]{P}):

\bigskip

\noindent \textbf{Lemma 3.6.} \textsl{The map }$\kappa $ \textsl{induces a
ring isomorphism (compare with (1.1))}

\begin{enumerate}
\item[(3.18)] $E_{3}^{\ast ,\ast }(U(n))=\Lambda (\xi _{1}^{\prime },\cdots
,\xi _{2n-1}^{\prime })\rightarrow H^{\ast }(U(n))=\Lambda (\xi _{1},\cdots
,\xi _{2n-1})$\textsl{.}\hfill $\square $
\end{enumerate}

\section{A refinement of the Gysin sequence (1.3)}

In general, the major difficulty to work with the Gysin sequence \cite[p.149]%
{MS} of a spherical fibration $c:E\rightarrow B$, such as (1.3), is to
evaluate the operator $\theta $ (e.g. Hirsch \cite{H}, Massey \cite{Ma}). It
requires, essentially, computing with the induced map $c^{\#}:C^{\ast
}(B)\rightarrow C^{\ast }(E)$ of $c$ at the co-chain level \cite[p.349-359]%
{Wh}. However, resorting to the Serre spectral sequence of the torus
fibrations $G\rightarrow G/T$ for $G=U(n)$ and $PU(n)$, we have concretly
interpreted the map $c^{\#}$ as a map of $d_{2}$-complexes

\begin{quote}
$c^{\ast }:$\textsl{\ }$E_{2}^{\ast ,\ast }(PU(n))\rightarrow E_{2}^{\ast
,\ast }(U(n))$ (see in Lemma 3.3).
\end{quote}

\noindent In this section we derive from $c^{\ast }$ an exact sequence
relating the groups $E_{3}^{\ast ,\ast }(U(n))$ and $E_{3}^{\ast ,\ast
}(PU(n))$, which serves also as a refinement of the Gysin sequence (1.3).
Moreover, extending the construction in Definition 3.5 we introduce the
primary $1$-forms $\rho _{3}^{\prime }$, $\cdots $, $\rho _{2n-1}^{\prime
}\in E_{3}^{\ast ,1}(PU(n))$ for the group $PU(n)$. They are applied in
Lemma 4.3 to derive preliminary versions of Theorems A and B.\hfill

With the product inherited from $E_{2}^{\ast ,\ast }(G)$ the third page $%
E_{3}^{\ast ,\ast }(G)$ is naturally a bi-graded ring \cite[P.668]{Wh}. In
term of i) of Lemma 3.3 the groups $E_{2}^{\ast ,r}(U(n))$ admits the
decomposition

\begin{quote}
$E_{2}^{\ast ,r}(U(n))=E_{2}^{\ast ,r}(PU(n))\oplus E_{2}^{\ast
,r-1}(PU(n))\otimes t_{0}$,
\end{quote}

\noindent so that the map $c^{\ast }$ fits into the short exact sequence of $%
d_{2}$-complexes

\begin{enumerate}
\item[(4.1)] $0\rightarrow E_{2}^{\ast ,r}(PU(n))\overset{c^{\ast }}{%
\rightarrow }E_{2}^{\ast ,r}(U(n))\overset{\theta }{\rightarrow }E_{2}^{\ast
,r-1}(PU(n))\rightarrow 0$,
\end{enumerate}

\noindent where, if $z,z^{\prime }\in E_{2}^{\ast ,\ast }(PU(n))$, then

\begin{quote}
$c^{\ast }(z)=z\oplus 0$,\quad $\theta (z\oplus z^{\prime }\otimes
t_{0})=z^{\prime }$.
\end{quote}

\noindent Conversely, every element $x\in E_{2}^{\ast ,r}(U(n))$\ can be
uniquely expressed as

\begin{enumerate}
\item[(4.2)] $x=x_{1}\oplus \theta (x)\otimes t_{0}$, where $x_{1}=x-\theta
(x)\otimes t_{0}\in E_{2}^{\ast ,r}(PU(n))$.
\end{enumerate}

\bigskip

\noindent \textbf{Lemma 4.1.} \textsl{The induced map of }$c^{\ast }$\textsl{%
\ fits into} \textsl{the exact sequence}

\begin{enumerate}
\item[(4.3)] {\small {$\cdots \rightarrow E_{3}^{\ast ,r}(PU(n))\overset{%
c^{\ast }}{\rightarrow }E_{3}^{\ast ,r}(U(n))\overset{\overline{\theta }}{%
\rightarrow }E_{3}^{\ast ,r-1}(PU(n))\overset{\overline{\omega }\cdot }{%
\rightarrow }E_{3}^{\ast ,r+1}(PU(n))\overset{c^{\ast }}{\rightarrow }\cdots 
$} }
\end{enumerate}

\noindent \textsl{where, if }$z\in E_{3}^{\ast ,\ast }(PU(n))$ \textsl{and }$%
x\in E_{3}^{\ast ,\ast }(U(n))$\textsl{,} \textsl{then}

\begin{quote}
\textsl{i)} $\overline{\theta }(c^{\ast }(z)\cdot x)=z\cdot \overline{\theta 
}(x)$\textsl{;}

\textsl{ii) }$\overline{\theta }(x\cdot x^{\prime })=x_{1}\cdot \overline{%
\theta }(x^{\prime })+(-1)^{\deg x^{\prime }}\overline{\theta }(x)\cdot
x_{1}^{\prime }$;\textsl{\quad }

\textsl{iii)} $\overline{\omega }(z)=\overline{\omega }\cdot z$.
\end{quote}

\noindent \textbf{Proof.} (4.3) is the exact sequence associated to the
short exact sequence (4.1) of $d_{2}$-complexes.

To see i) consider two arbitrary elements

\begin{quote}
$x=x_{1}\oplus \theta (x)\otimes t_{0}\in E_{2}^{\ast ,r}(U(n))$ and $z\in
E_{2}^{\ast ,r}(PU(n))$.
\end{quote}

\noindent Then the expression (4.2) of the product $c^{\ast }(z)\cdot x$
turns to be

\begin{quote}
$c^{\ast }(z)\cdot x=(z\oplus 0)(x_{1}\oplus \theta (x)\otimes t_{0})=z\cdot
x_{1}\oplus z\cdot \theta (x)\otimes t_{0}.$
\end{quote}

\noindent Relation i) is verified by the definition of $\theta $.

For ii) take two arbitrary elements of $E_{2}^{\ast ,\ast }(U(n))$

\begin{quote}
$x=x_{1}\oplus \theta (x)\otimes t_{0}$, $x^{\prime }=x_{1}^{\prime }\oplus
\theta (x^{\prime })\otimes t_{0}$.
\end{quote}

\noindent Then, by the multiplicative properties in $E_{2}^{\ast ,\ast
}(U(n))$

\begin{quote}
$\left( 1\otimes t_{0}\right) ^{2}=0$ and $a\cdot b=(-1)^{\deg a\deg
b}b\cdot a$,
\end{quote}

\noindent the expression (4.2) of the product $x\cdot x^{\prime }$ is

\begin{quote}
$x\cdot x^{\prime }=x_{1}\cdot x_{1}^{\prime }\oplus (x_{1}\cdot \theta
(x^{\prime })+(-1)^{\deg x^{\prime }}\theta (x)\cdot x_{1}^{\prime })\otimes
t_{0}$.
\end{quote}

\noindent We obtain ii) again by the definition of $\theta $.

Finally, for a $d_{2}^{\prime }$-cocycle $z\in E_{2}^{\ast ,r-1}(PU(n))$
consider the diagram chasing in the short exact sequence (4.1) in $d_{2}$%
-complexes:

\begin{quote}
$%
\begin{array}{ccccc}
&  & 0\oplus z\otimes t_{0} & \overset{\theta }{\rightarrow } & z \\ 
&  & d_{2}\downarrow \quad &  & d_{2}^{\prime }\downarrow \quad \\ 
z\cdot x_{1} & \overset{c^{\ast }}{\rightarrow } & z\cdot x_{1}\oplus 0 &  & 
0%
\end{array}%
$.
\end{quote}

\noindent With $\overline{\omega }=[x_{1}\otimes 1]$ (see in Lemma 3.3) it
says that the connecting homomorphism $\beta $ in the exact sequence (4.3)
is $\beta ([z])=[z]\cdot \overline{\omega }$. This shows iii), completing
the proof of the lemma.\hfill $\square $

\bigskip

In the exact sequence (4.3), the groups $E_{3}^{\ast ,r}(U(n))$ are known
(i.e. Lemma 3.6). To derive the unknown groups $E_{3}^{\ast ,r}(PU(n))$ from
(4.3) we construct, in the same spirit of Definition 3.5, a set of primary
1-forms $\rho _{3}^{\prime }$, $\cdots $, $\rho _{2n-1}^{\prime }\in
E_{3}^{\ast ,1}(PU(n))$ for the group $PU(n)$. Consider again the lift $%
\widehat{e}_{k}$ of the symmetric polynomial $e_{k}\in \mathbb{Z}%
[x_{1},\cdots ,x_{n}]$ with respect to $d_{2}$

\begin{quote}
$\widehat{e}_{k}=e_{k}^{(1)}\otimes t_{0}+e_{k}^{(2)}\otimes t_{1}+\cdots
+e_{k}^{(n)}\otimes t_{n-1}$.
\end{quote}

\noindent We have in the notation of (4.2) that

\begin{quote}
$\widehat{e}_{k}=\widehat{a}_{k}\oplus \theta (\widehat{e}_{k})\otimes
t_{0}\in E_{2}^{\ast ,1}(U(n))=E_{2}^{\ast ,1}(PU(n))\oplus E_{2}^{\ast
,0}(PU(n))\otimes t_{0}$,
\end{quote}

\noindent where

\begin{enumerate}
\item[(4.4)] $\widehat{a}_{k}:=e_{k}^{(2)}\otimes t_{1}+\cdots
+e_{k}^{(n)}\otimes t_{n-1}\in E_{2}^{2k-1,1}(PU(n))$;

$\theta (\widehat{e}_{k})=e_{k}^{(1)}=\binom{n}{k}x_{1}^{k-1}\in
E_{2}^{2(k-1),0}(PU(n))$ (by (3.14)).
\end{enumerate}

\noindent From $d_{2}(\widehat{e}_{k})=d_{2}^{\prime }(\widehat{a}_{k})+%
\binom{n}{k}x_{1}^{k}=0$ by ii) of Lemma 3.3 and (3.16) we get

\begin{enumerate}
\item[(4.5)] $d_{2}^{\prime }(\widehat{a}_{k})=-\binom{n}{k}x_{1}^{k}\in
E_{2}^{2k,0}(PU(n))$.
\end{enumerate}

\noindent In addition, for each $1\leq k\leq n$, the relation $b_{n,k}%
\overline{\omega }^{k}=0$ in (3.9) implies that there exists an element $%
\widehat{\delta }_{n,k}\in E_{2}^{2(k-1),1}(PU(n))$\textsl{\ }so that

\begin{enumerate}
\item[(4.6)] $d_{2}^{\prime }(\widehat{\delta }_{n,k})=b_{n,k}\cdot
x_{1}^{k} $.
\end{enumerate}

Suppose that the integer $n$ has the prime factorization $%
p_{1}^{r_{1}}\cdots p_{t}^{r_{t}}$, and recall from Section \S 1 the
partition of the set of integers $\left\{ 2,\cdots ,n\right\} $:

\begin{quote}
$\left\{ 2,\cdots ,n\right\} =Q_{0}(n)\underset{1\leq i\leq t}{\amalg }%
Q_{p_{i}}(n)$, where $Q_{p_{i}}(n):=\{p_{i},\cdots ,p_{i}^{r_{i}}\}$.
\end{quote}

\bigskip

\noindent \textbf{Definition 4.2.} If $k\in Q_{0}(n)$ then $%
b_{n,k}=b_{n,k-1} $ by i) of Lemma 2.1. In particular, the element

\begin{quote}
$\widehat{h}_{k}:=\widehat{a}_{k}+\frac{\binom{n}{k}}{b_{n,k}}\cdot
x_{1}\cdot \widehat{\delta }_{n,k-1}\in E_{2}^{2(k-1),1}(PU(n))$
\end{quote}

\noindent is a $d_{2}^{\prime }$-cocycle by (4.5) and (4.6). We obtain
therefore the cohomology class

\begin{enumerate}
\item[(4.7)] $\rho _{2k-1}^{\prime }:=[\widehat{h}_{k}]\in
E_{3}^{2(k-1),1}(PU(n))$.
\end{enumerate}

If $k=p^{s}\in Q_{p}(n)$ with $p\in \{p_{1},\cdots ,p_{t}\}$, consider the
element of $E_{2}^{2(k-1),1}(PU(n))$

\begin{quote}
$\widehat{h}_{k}:=p\cdot \widehat{a}_{p^{s}}-\frac{n}{p^{s-1}}\underset{%
1\leq t\leq p^{s}-p^{s-1}}{\Sigma }(-1)^{t-1}x_{1}^{t}\cdot \widehat{a}%
_{p^{s}-t}-x_{1}^{p^{s}-p^{s-1}}\cdot \widehat{a}_{p^{s-1}}$.
\end{quote}

\noindent Since $d_{2}^{\prime }(\widehat{h}_{k})=0$ by (4.5) and (2.2), we
obtain the cohomology class

\begin{enumerate}
\item[(4.8)] $\rho _{2k-1}^{\prime }:=[\widehat{h}_{k}]\in
E_{3}^{2(r-1),1}(PU(n))$.
\end{enumerate}

\noindent The elements $\rho _{3}^{\prime },\cdots ,\rho _{2n-1}^{\prime }$
obtained in (4.7) and (4.8) are called \textsl{the} \textsl{primary }$1$%
\textsl{-forms} \textsl{of} the group $PU(n)$.\hfill $\square $

\bigskip

Granted with the primary $1$-forms $\xi _{1}^{\prime },\cdots ,\xi
_{2n-1}^{\prime }$ for $G=U(n)$ and $\rho _{3}^{\prime },\cdots ,\rho
_{2n-1}^{\prime }$ for $G=$ $PU(n)$, respectively, we establish below a
preliminary version of Theorems A and B by computation in the much simpler
rings $E_{2}^{\ast ,\ast }(G)$ (i.e. Lemma 3.3). For a multi-index $I=$%
\textsl{\ }$\{i_{1},\cdots ,i_{k}\}\subseteq $\textsl{\ }$\{1,\cdots ,n\}$%
\textsl{\ }with\textsl{\ }$k\geq 2$ we put

\begin{quote}
$I^{e}:=${\small \ }$\{i_{1},\cdots ,i_{k-1},\widehat{i_{k}}\}$, $\xi
_{I}^{\prime }:=\underset{j\in I}{\cup }${\small \ }$\xi _{2j-1}^{\prime }$.
\end{quote}

\noindent \textbf{Lemma 4.3.} \textsl{In the exact sequence (4.3) the }$1$%
\textsl{-forms }$\xi _{2k-1}^{\prime }$ \textsl{and} $\rho _{2k-1}^{\prime }$%
\textsl{\ satisfy the following relations:}

\begin{enumerate}
\item[i)] $c^{\ast }(\rho _{2k-1}^{\prime })=$ $\xi _{2k-1}^{\prime }$ 
\textsl{if }$k\in Q_{0}(n)$\textsl{;}

\item[ii)] $c^{\ast }(\rho _{2k-1}^{\prime })=p\cdot \xi _{2k-1}^{\prime }$ 
\textsl{if }$k\in Q_{p}(n)$\textsl{;}

\item[iii)] $\overline{\theta }(\xi _{2k-1}^{\prime })=\binom{n}{k}\overline{%
\omega }^{k-1}$\textsl{, }$1\leq k\leq n$\textsl{;}

\item[iv)] $\overline{\theta }(\xi _{I}^{\prime })=-\overline{\theta }(\xi
_{I^{e}}^{\prime })\cdot \rho _{2i_{k}-1}$ \textsl{if} $i_{k}\in Q_{0}(n)$%
\textsl{;}

\item[v)] $\overline{\theta }(\xi _{I}^{\prime })=-\frac{1}{p}(\overline{%
\theta }(\xi _{I^{e}}^{\prime })\cdot \rho _{2p^{s}-1}-\overline{\omega }%
^{p^{s}-p^{s-1}}\cdot \overline{\theta }(\xi _{I^{e}}^{\prime }\cdot \xi
_{2p^{s-1}-1}^{\prime }))$ \textsl{if} $i_{k}=p^{s}\in Q_{p}(n)$\textsl{.}
\end{enumerate}

\noindent \textbf{Proof.} For $k\in Q_{0}(n)$ define $b_{k}:=\frac{\binom{n}{%
k}}{b_{n,k}}\widehat{\delta }_{n,k-1}\otimes t_{0}\in E_{2}^{2(k-2),2}(U(n))$%
. From

\begin{quote}
$d_{2}(b_{k})=\frac{\binom{n}{k}b_{n,k-1}}{b_{n,k}}x_{1}^{k-1}\otimes t_{0}-%
\frac{\binom{n}{k}}{b_{n,k}}x_{1}\widehat{\delta }_{n,k-1}\otimes 1$
\end{quote}

\noindent and $b_{n,k-1}=b_{n,k}$ by i) of Lemma 2.1, one finds that

\begin{quote}
$\widehat{e}_{k}-c^{\ast }(\widehat{h}_{k})=d_{2}(b_{k})$.
\end{quote}

\noindent One obtains i) from $\xi _{2k-1}^{\prime }=\left[ \widehat{e}_{k}%
\right] $, $\rho _{2k-1}^{\prime }=[\widehat{h}_{k}]$ and $[d_{2}(b_{k})]=0$.

Similarly, for the case $k=p^{s}\in Q_{p}(n)$ with $p\in \{p_{1},\cdots
,p_{t}\}$ we take

\begin{quote}
$b_{k}:=(\frac{n}{p^{s-1}}\underset{1\leq t\leq p^{s}-p^{s-1}}{\Sigma }%
(-1)^{t-1}x_{1}^{t-1}\cdot \widehat{a}_{p^{s}-t}+(-x_{1})^{p^{s}-p^{s-1}}%
\cdot \widehat{a}_{p^{s-1}})\otimes t_{0}$.
\end{quote}

\noindent Since $p\cdot \widehat{e}_{k}-c^{\ast }(\widehat{h}%
_{k})=d_{2}(b_{k})$ by (4.5) and (4.8), we obtain ii).

Relation iii) follows directly from $\theta (\widehat{e}_{k})=\binom{n}{k}%
x_{1}^{k-1}$ by (4.4). To see the relations iv) and v) we assume below that $%
I=$\textsl{\ }$\{i_{1},\cdots ,i_{k}\}$ with $k\geq 2$.

If $i_{k}\in Q_{0}(n)$, we get by i) of Lemma 4.1 and $c^{\ast }(\rho
_{2i_{k}-1}^{\prime })=\xi _{2i_{k}-1}^{\prime }$ that

\begin{quote}
$\overline{\theta }(\xi _{I}^{\prime })=\overline{\theta }(\xi
_{I^{e}}^{\prime }\cdot c^{\ast }(\rho _{2i_{k}-1}^{\prime }))=-\overline{%
\theta }(\xi _{I^{e}}^{\prime })\cdot \rho _{2i_{k}-1}$,
\end{quote}

\noindent showing the relation iv).

Finally, for the case $i_{k}=p^{s}\in Q_{p}(n)$ we set

\begin{quote}
$\widehat{e}_{I}:=\underset{i\in I}{\Pi }\widehat{e}_{i}\in E_{2}^{\ast
,k}(U(n))$,$\quad \widehat{a}_{I}:=\underset{i\in I}{\Pi }\widehat{a}_{i}\in
E_{2}^{\ast ,k}(PU(n))$.
\end{quote}

\noindent In particular, we have $\widehat{e}_{I}=\widehat{a}_{I}\oplus
\theta (\widehat{e}_{I})\otimes t_{0}$ by (4.2). From

\begin{quote}
$d_{2}(\widehat{e}_{I})=0$, $d_{2}^{\prime }(\theta (\widehat{e}_{I}))=0$
and $d_{2}(1\otimes t_{0})=x_{1}$,
\end{quote}

\noindent we get

\begin{enumerate}
\item[(4.9)] $d_{2}^{\prime }(\widehat{a}_{I})=(-1)^{k}x_{1}\cdot \theta (%
\widehat{e}_{I})$.
\end{enumerate}

\noindent On the other hand, according to ii) of Lemma 4.1 we have

\begin{quote}
$\theta (\widehat{e}_{I})=\theta (\widehat{e}_{I^{e}}\widehat{e}_{p^{s}})=%
\widehat{a}_{I^{e}}\cdot \binom{n}{p^{s}}x_{1}^{p^{s}-1}-\theta (\widehat{e}%
_{I^{e}})\cdot \widehat{a}_{p^{s}}$,

$\theta (\widehat{e}_{I^{e}}\widehat{e}_{p^{s-1}})=\widehat{a}_{I^{e}}\cdot 
\binom{n}{p^{s-1}}x_{1}^{p^{s-1}-1}-\theta (\widehat{e}_{I^{e}})\cdot 
\widehat{a}_{p^{s-1}}$.
\end{quote}

\noindent Using these formulae, a direct but tedious computation shows that

\begin{quote}
$\theta (\widehat{e}_{I})+\frac{1}{p}\theta (\widehat{e}_{I^{e}})\cdot 
\widehat{h}_{p^{s}}-\frac{1}{p}x_{1}^{p^{s}-p^{s-1}}\cdot \theta (\widehat{e}%
_{I^{e}}\cdot \widehat{e}_{p^{s-1}})$

$=d_{2}^{\prime }(b_{I})$ (by (2.2) and (4.9)),
\end{quote}

\noindent where

\begin{quote}
$b_{I}:=\frac{n}{p^{s}}\underset{1\leq t\leq p^{s}-p^{s-1}}{\Sigma }%
(-1)^{t-1}x_{1}^{t-1}\cdot \widehat{a}_{I^{e}}\cdot \widehat{a}_{p^{s}-t}\in
E_{2}^{\ast ,k}(PU(n))$.
\end{quote}

\noindent This completes the proof of the relation v).\hfill \hfill $\square 
$

\section{Proofs of Theorem A and B}

Applying the map $\kappa $ in (3.12) to the primary $1$-forms $\rho
_{3}^{\prime },\cdots ,\rho _{2n-1}^{\prime }\in E_{3}^{\ast ,1}(PU(n))$
introduced in Definition 4.2, we define the desired cohomology classes $\rho
_{2k-1}$ stated in Theorem A by

\begin{enumerate}
\item[(5.1)] $\rho _{2k-1}:=\kappa (\rho _{2k-1}^{\prime })\in
H^{2k-1}(PU(n))$, $2\leq k\leq n$.
\end{enumerate}

\noindent Since the map $\kappa $ is natural with respect to the bundle map $%
c$, we have the commutative diagram

\begin{quote}
\begin{tabular}{lll}
$E_{3}^{\ast ,1}(PU(n))$ & $\overset{c^{\ast }}{\rightarrow }$ & $%
E_{3}^{\ast ,1}(U(n))$ \\ 
$\kappa \downarrow $ &  & $\kappa \downarrow $ \\ 
$H^{\ast }(PU(n))$ & $\overset{c^{\ast }}{\rightarrow }$ & $H^{\ast }(U(n))$.%
\end{tabular}
\end{quote}

\noindent It implies, by the relations i) and ii) of Lemma 4.3, that

\bigskip

\noindent \textbf{Corollary 5.1.} \textsl{The map }$c^{\ast }:H^{\ast
}(PU(n))\rightarrow H^{\ast }(U(n))$\textsl{\ satisfies that}

\begin{enumerate}
\item[(5.2)] $c^{\ast }(\rho _{2k-1})=\left\{ 
\begin{array}{c}
\xi _{2k-1}\text{ \textsl{if} }k\in Q_{0}(n)\text{\textsl{;}} \\ 
p_{i}\cdot \xi _{2k-1}\text{ \textsl{if} }k\in Q_{p_{i}}(n)\text{\textsl{.}}%
\end{array}%
\right. $
\end{enumerate}

\noindent \textsl{In particular,} $c^{\ast }(\rho _{3}\cdots \rho
_{2n-1})=n\cdot \xi _{3}\cdots \xi _{2n-1}$\textsl{.}\hfill \hfill $\square $

\bigskip

Let $\left\langle \omega \right\rangle \subset $ $H^{\ast }(PU(n))$ be the
ideal generated by the Euler class $\omega $, and denote by $H^{\ast
}(PU(n))_{\left\langle \omega \right\rangle }$ the quotient ring $H^{\ast
}(PU(n))/\left\langle \omega \right\rangle $. Then the Gysin sequence (1.3)
can be summarized into the exact sequences with four terms

\begin{enumerate}
\item[(5.3)] $0\rightarrow H^{\ast }(PU(n))_{\left\langle \omega
\right\rangle }\overset{c^{\ast }}{\rightarrow }H^{\ast }(U(n))\overset{%
\theta }{\rightarrow }H^{\ast }(PU(n))\overset{\omega \cdot }{\rightarrow }%
\left\langle \omega \right\rangle \rightarrow 0$.
\end{enumerate}

\noindent Since the ring $H^{\ast }(U(n))$ is torsion free, while the ring
map $c^{\ast }$ in (5.3) injects, the quotient ring $H^{\ast
}(PU(n))_{\left\langle \omega \right\rangle }$ must be torsion free.
Therefore, the short exact sequence

\begin{quote}
$0\rightarrow \left\langle \omega \right\rangle \rightarrow H^{\ast
}(PU(n))\rightarrow H^{\ast }(PU(n))_{\left\langle \omega \right\rangle
}\rightarrow 0$
\end{quote}

\noindent is splittable to yield the decomposition

\begin{enumerate}
\item[(5.4)] $H^{\ast }(PU(n))=H^{\ast }(PU(n))_{\left\langle \omega
\right\rangle }\oplus \left\langle \omega \right\rangle $.
\end{enumerate}

\noindent In addition, as the order of the Euler class $\omega $ is
precisely $n$ by Lemma 3.4,\textbf{\ }we can conclude by (5.4) that the two
summands $H^{\ast }(PU(n))_{\left\langle \omega \right\rangle }$ and $%
\left\langle \omega \right\rangle $ are the free part and the torsion part
of the cohomology ring $H^{\ast }(PU(n))$, respectively.

\bigskip

\noindent \textbf{Lemma 5.2. }\textsl{With the notation above} \textsl{one
has}

\textsl{i) }$\rho _{2k-1}^{2}=0$\textsl{;}

\textsl{ii)} \textsl{The subring of }$H^{\ast }(PU(n))$\textsl{\ generated
by }$\rho _{3},\cdots ,\rho _{2n-1}$\textsl{\ is isomorphic to the exterior
ring }$\Lambda (\rho _{3},\cdots ,\rho _{2n-1})$\textsl{, which represents
also} \textsl{the free part} $H^{\ast }(PU(n))_{\left\langle \omega
\right\rangle }$\textsl{\ of }$H^{\ast }(PU(n))$\textsl{.}

\bigskip

\noindent \textbf{Proof. }Let $\beta _{2}:$ $H^{r}(PU(n);\mathbb{Z}%
_{2})\rightarrow H^{r+1}(PU(n))$ be the $\func{mod}2$ Bockstein
homomorphism, and assume that $n=2^{m}(2b+1)$. It follows from the formula
(5.7) below that

\begin{quote}
$\func{Im}\pi ^{\prime \ast }\cap \func{Im}\beta _{2}=\{0,2^{m-t-1}\omega
^{2^{t}}\mid 0\leq t\leq m-1\}$.
\end{quote}

\noindent (see also \cite[Lemma 6.7]{BB} or \cite[Proposition 4.2]{P2}). On
the other hand, it is known \cite[Lemma 2.2]{D2} that, for any compact Lie
group $G$ and any $\rho \in E_{3}^{\ast ,1}(G)$,

\begin{quote}
$\kappa (\rho )^{2}\in \func{Im}\pi ^{\ast }\cap \func{Im}\beta _{2}$.
\end{quote}

\noindent We obtain $\rho _{2k-1}^{2}=0$ from $\deg \omega =2$ but $\deg
\rho _{2k-1}^{2}=2(2k-1)$.

The statements in ii) comes directly from the relation

\begin{quote}
$c^{\ast }(\rho _{3}\cdots \rho _{2n-1})=n\cdot \xi _{3}\cdots \xi _{2n-1}$
\end{quote}

\noindent by Corollary 5.1, which have been shown in \cite[Lemma 2.2]{DL}%
.\hfill $\square $

\bigskip

\noindent \textbf{Proof of Theorem A.} The relation (1.4) in Theorem A has
been shown by (5.2). It remains to show isomorphism (1.5).

In view of (5.4) and by ii) of Lemma 5.2, every element $z\in H^{\ast
}(PU(n))$ can be expressed as

\begin{quote}
$z=a_{0}\oplus \omega \cup z_{1}$
\end{quote}

\noindent for some $a_{0}\in \Lambda (\rho _{3},\cdots ,\rho _{2n-1})$ and $%
z_{1}\in H^{\ast }(PU(n))$. By repetition and for the degree reason, we get
the expansion

\begin{quote}
$z=a_{0}+\omega \cdot a_{1}+\omega ^{2}\cdot a_{2}+\cdots $ $+$ $\omega
^{k}\cdot a_{k}$, $a_{i}\in \Lambda (\rho _{3},\cdots ,\rho _{2n-1})$,
\end{quote}

\noindent indicating that the inclusion $\{\omega ,\rho _{2k-1}\}\subset
H^{\ast }(PU(n))$ extends to an epimorphism

\begin{quote}
$h:\mathbb{Z}[\omega ]\otimes \Lambda (\rho _{3},\cdots ,\rho
_{2n-1})\rightarrow H^{\ast }(PU(n))$.
\end{quote}

\noindent It follows from $\ker (\omega \cdot )=\func{Im}\theta $ (by the
exactness of the sequence (5.3)) that

\begin{quote}
$\ker h=\langle \omega \otimes \theta (\xi _{I}),I\subseteq \{1,\cdots
,n\}\rangle $,
\end{quote}

\noindent completing the proof of Theorem A.\hfill $\square $

\bigskip

Based on the exact sequence (4.3), together with properties i) and ii) of
Lemma 4.4, the same argument proving Theorem A shows that

\bigskip

\noindent \textbf{Lemma 5.3. }\textsl{Let}\textbf{\ }$\overline{\theta }$ 
\textsl{be the map} $E_{3}^{\ast ,\ast }(U(n))\rightarrow E_{3}^{\ast ,\ast
}(PU(n))$ \textsl{in (4.3). Then}

\begin{enumerate}
\item[(5.5)] $E_{3}^{\ast ,\ast }(PU(n))=\mathbb{Z}[\overline{\omega }%
]\otimes \Lambda (\rho _{3}^{\prime },\cdots ,\rho _{2n-1}^{\prime
})/\langle \overline{\omega }\otimes \overline{\theta }(\xi _{I}^{\prime
}),I\subseteq \{1,\cdots ,n\}\rangle $\textsl{.}\hfill $\square $
\end{enumerate}

\bigskip

Comparing formula (1.5) with (5.5), we obtain the following outcome, which
is consistent with the Borel-Pittie isomorphism (3.18) for the group $G=U(n)$%
.

\bigskip

\noindent \textbf{Corollary 5.4.} \textsl{The map }$K:E_{3}^{\ast ,\ast
}(PU(n))\rightarrow H^{\ast }(PU(n))$\textsl{\ defined by }$\pi ^{\prime
\ast }(\overline{\omega })=\omega $\textsl{\ and} $\kappa (\rho
_{2k-1}^{\prime })=\rho _{2k-1}$ \textsl{is an isomorphism.}\hfill $\square $

\bigskip

\noindent \textbf{Proof of Theorem B.} By the naturalities of the maps $\pi
^{\ast }$ and $\kappa $ with respect to the bundle map $c$, we have the
commutative diagrams in the maps $L^{^{\prime }}$, $K^{^{\prime }}$ between
the exact sequences (4.3) and (1.3)

\begin{center}
{\small {%
\begin{tabular}{lllllll}
$\cdots \rightarrow $ & $E^{p,q}(PU(n))$ & $\overset{c^{\ast }}{\rightarrow }
$ & $E^{p,q}(U(n))$ & $\overset{\overline{\theta }}{\rightarrow }$ & $%
E^{p,q-1}(PU(n))$ & $\overset{\overline{\omega }\cdot }{\rightarrow }\cdots $
\\ 
& $K^{^{\prime }}\downarrow $ &  & $L^{^{\prime }}\downarrow $ &  & $%
K^{^{\prime }}\downarrow $ &  \\ 
$\cdots \rightarrow $ & $H^{p+q}(PU(n))$ & $\overset{c^{\ast }}{\rightarrow }
$ & $H^{p+q}(U(n))$ & $\overset{\theta }{\rightarrow }$ & $H^{p+q-1}(PU(n))$
& $\overset{\omega \cup }{\rightarrow }\cdots $%
\end{tabular}%
}, }
\end{center}

\noindent where $K^{^{\prime }}$ is the restriction of the isomorphism $K$
to the subgroup $E_{3}^{p,q}(PU(n))$ $\subset E_{3}^{\ast ,\ast }(PU(n))$,
and $L^{^{\prime }}$ is the restriction of the isomorphism (3.18) to the
subgroup $E^{p,q}(U(n))$. It allows us to translate the properties iii),
iv), v) of Lemma 4.3 to the desired relations a), b), c) stated in Theorem
B.\hfill $\square $

\bigskip

Let $G$ be a compact connected Lie group. In \cite{G} Grothendieck
introduced the \textsl{Chow ring} of the reductive algebraic group $G^{c}$
corresponding to $G$ by

\begin{quote}
$\mathcal{A}^{\ast }(G^{c}):=\func{Im}\{\pi ^{\ast }:H^{\ast
}(G/T)\rightarrow H^{\ast }(G)\}$.
\end{quote}

\noindent For $G=PU(n)$ let $J_{n}(\omega )\subset H^{\ast }(PU(n))$ be the
subring generated by the Euler class $\omega \in H^{2}(PU(n))$. By Corollary
5.4 and Lemma 3.4, the map

\begin{quote}
$\pi ^{^{\prime }\ast }:E_{3}^{\ast ,0}(PU(n))\rightarrow H^{\ast }(PU(n))$
\end{quote}

\noindent in (3.10) induces an isomorphism $E_{3}^{\ast ,0}(PU(n))\cong
J_{n}(\omega )$. We get by formula (3.9) that

\bigskip

\noindent \textbf{Corollary 5.5. }\textsl{The Chow ring of the algebraic
group }$PU(n)^{c}$\textsl{\ is isomorphic to}

\begin{enumerate}
\item[(5.6)] $\func{Im}\pi ^{^{\prime }\ast }=J_{n}(\omega )=\frac{\mathbb{Z}%
[{\small \omega }]}{\left\langle b_{n,k}\omega ^{k},1\leq k\leq
n\right\rangle }$\textsl{.}
\end{enumerate}

\textsl{In particular,} \textsl{if }$n>2$ \textsl{has} \textsl{the prime
factorization }$p_{1}^{r_{1}}\cdots p_{t}^{r_{t}}$\textsl{, then one has the
splitting}

\begin{enumerate}
\item[(5.7)] $J_{n}(\omega )=\mathbb{Z\oplus (}\underset{i\in \{1,\cdots ,t\}%
}{\mathbb{\oplus }}$ $\frac{\mathbb{Z}[\omega ]^{+}}{\left\langle
p_{i}^{r_{i}}\omega ,p_{i}^{r_{i}-1}\omega ^{p_{i}},p_{i}^{r_{i}-2}\omega
^{p_{i}^{2}},\cdots ,\omega ^{p_{i}^{r_{i}}}\right\rangle })$,
\end{enumerate}

\noindent \textsl{where }$\mathbb{Z}[\omega ]^{+}$\textsl{\ is the subring
of }$\mathbb{Z}[\omega ]$\textsl{\ consisting of polynomials in }$\omega $%
\textsl{\ without constant terms.}

\bigskip

\noindent \textbf{Proof.} It suffices to show (5.7). As is customary, any
finitely generated graded group $B$ admits the decomposition

\begin{quote}
$B=\mathcal{F}(B)\oplus (\underset{p}{\oplus }\sigma _{p}(B))$,
\end{quote}

\noindent where $\mathcal{F}(B):=B/TorB$ is \textsl{the free part} of $B$,
the sum is over all prime $p\geq 2$, and $\sigma _{p}(B)$ is \textsl{the} $p$%
\textsl{-primary component} of $B$ defined by

\begin{quote}
$\sigma _{p}(B):=\{x\in B,$ $p^{r}x=0,$ $r\geq 1\}$.
\end{quote}

\noindent For the case $B=J_{n}(\omega )$ one gets by (5.6) and ii) of Lemma
2.1 that if $p\notin \{p_{1},\cdots ,p_{t}\}$, then $\sigma
_{p}(J_{n}(\omega ))=0$, and if $p=p_{i}\in \{p_{1},\cdots ,p_{t}\}$, then

\begin{quote}
$\sigma _{p}(J_{n}(\omega ))=\frac{\mathbb{Z}[\omega ]^{+}}{\left\langle
p_{i}^{r_{i}}\omega ,p_{i}^{r_{i}-1}\omega ^{p_{i}},p_{i}^{r_{i}-2}\omega
^{p_{i}^{2}},\cdots ,\omega ^{p_{i}^{r_{i}}}\right\rangle }$.
\end{quote}

\noindent That is, formula (5.7) follows from (5.6).\hfill $\square $

\section{Main result and examples}

For any finite $CW$-complex $X$ the cohomology $H^{\ast }(X)$ has the
decomposition

\begin{enumerate}
\item[(6.1)] $H^{\ast }(X)=\mathcal{F}(X)\oplus (\underset{p}{\oplus }\sigma
_{p}(X))$,
\end{enumerate}

\noindent where $\mathcal{F}(X)$ and $\sigma _{p}(X)$ are the free part and
the $p$-primary component of $H^{\ast }(X)$, respectively (see in the proof
of Corollary 5.5). In this section, we derive such an expression of the ring 
$H^{\ast }(PU(n))$, by describing its basic components $\mathcal{F}(PU(n))$
and $\sigma _{p}(PU(n))$ explicitly.

Suppose, therefore, that the prime factorization of $n$ is $%
p_{1}^{r_{1}}\cdots p_{t}^{r_{t}}$, and that $p=p_{i}\in \{p_{1},\cdots
,p_{t}\}$. A multi-index $I\subseteq \{1,\cdots ,n\}$ is called $p$\textsl{%
-monotone} if

\begin{quote}
$I\subseteq Q_{p}^{+}(n):=\{1,p,\cdots ,p^{r}\}$, where $r=r_{i}$.
\end{quote}

\noindent For a $p$-monotone sequence $I=\left\{ p^{i_{1}},\cdots
,p^{i_{k}}\right\} $ the two \textsl{derived sequences}

\begin{quote}
$I^{e}:=\left\{ p^{i_{1}},\cdots ,p^{i_{k-1}}\right\} $ and $I^{\partial
}:=\left\{ p^{i_{1}},\cdots ,p^{i_{k-1}},p^{i_{k}-1}\right\} $
\end{quote}

\noindent are clearly $p$-monotone, while the formula c) in Theorem B turns
to be

\begin{quote}
$\theta (\xi _{I})=-\frac{1}{p}(\theta (\xi _{I^{e}})\cup \rho
_{2p^{i_{k}}-1}+\omega ^{p^{i_{k}}-p^{i_{k}-1}}\cup \theta (\xi
_{I^{\partial }}))$.
\end{quote}

\noindent\ Repeatedly applying this to reduce the term $\theta (\xi
_{I^{\partial }})$ on the right side to get

\begin{enumerate}
\item[(6.2)] $\theta (\xi _{I})=-\frac{1}{p}\theta (\xi _{I^{e}})\cup (%
\underset{0\leq t\leq i_{k}-i_{k-1}-1}{\Sigma }\frac{1}{p^{t}}\omega
^{p^{i_{k}}-p^{i_{k}-t}}\cup \rho _{2p^{i_{k}-t}-1})$,
\end{enumerate}

\noindent where the sum must end at $t=i_{k}-i_{k-1}-1$, because the
relation $\xi _{2i-1}^{2}=0$ on $H^{\ast }(U(n))$ implies $\xi _{I^{\partial
}}=0$ whenever $i_{k}=i_{k-1}+1$. To see the implication of (6.2) we
introduce the set

\begin{quote}
$S(I):=\{J=\left\{ p^{j_{2}},\cdots ,p^{j_{k}}\right\} \subseteq
Q_{p}^{+}(n);i_{s}\geq j_{s}>i_{s-1}\}$,
\end{quote}

\noindent and define the functions $\varepsilon _{I}$, $\kappa
_{I}:S(I)\rightarrow \mathbb{Z}$ by

\begin{quote}
$\varepsilon _{I}(J):=(i_{k}-j_{k})+\cdots +(i_{2}-j_{2})$ and

$\kappa _{I}(J):=(p^{i_{k}}-p^{j_{k}})+\cdots
+(p^{i_{2}}-p^{j_{2}})+p^{i_{1}}-1$, $J\in S(I)$,
\end{quote}

\noindent respectively. Then, taking (6.2) as a new recurrence to reduce the
length of $I$, and applying the formula $\theta (\xi _{2k-1})=\binom{n}{k}%
\omega ^{k-1}$ by a) of Theorem B to evaluate the factor $\theta (\xi
_{2p^{i_{1}}-1})$ at the final step, we obtain the following formula, that
expresses $\theta (\xi _{I})$ as an explicit element of the ring $\mathbb{Z}%
[\omega ]\otimes \Lambda (\rho _{3},\cdots ,\rho _{2n-1})$.

\bigskip

\noindent \textbf{Lemma 6.1. }\textsl{If }$I=\left\{ p^{i_{1}},\cdots
,p^{i_{k}}\right\} $\textsl{\ }$\subseteq Q_{p}^{+}(n)$ \textsl{is} $p$%
\textsl{-monotone, then}

\begin{enumerate}
\item[(6.3)] $\theta (\xi _{I})=\frac{(-1)^{k-1}}{p^{r-i_{1}}}\binom{n}{%
p^{i_{1}}}\underset{J\in S(I)}{\Sigma }p^{r-i_{1}-\varepsilon _{I}(J)}\cdot
\omega ^{\kappa _{I}(J)}\cup \rho _{J}$\textsl{,}
\end{enumerate}

\noindent \textsl{where }$p^{r-i_{1}}$\textsl{\ divides }$\binom{n}{p^{i_{1}}%
}$\textsl{\ by ii) of Lemma 2.1.}\hfill $\square $

\bigskip

Observe that, for an arbitrary subsequence $I\subseteq \{1,\cdots ,n\}$ and $%
1\leq i\leq t$, the intersection $I_{i}:=Q_{p_{i}}^{+}(n)\cap I$ is either
empty, or is $p_{i}$-monotone. The following result allows one to express $%
\theta (\xi _{I})$ in term of $\theta (\xi _{I_{i}})$, whose formula is
available by (6.3).

\bigskip

\noindent \textbf{Lemma 6.2.} \textsl{For any multi--index }$I\subseteq
\{1,\cdots ,n\}$ \textsl{we have}

\begin{enumerate}
\item[(6.4)] $\theta (\xi _{I})=a_{I}^{1}\cup \theta (\xi _{I_{1}})+\cdots
+a_{I}^{t}\cup \theta (\xi _{I_{t}})$\textsl{\ with} $a_{I}^{i}\in \Lambda
(\rho _{3},\cdots ,\rho _{2n-1})$\textsl{.}
\end{enumerate}

\textsl{In particular, the ring epimorphism in the proof of Theorem A}

\begin{quote}
$h:\mathbb{Z}[\omega ]\otimes \Lambda (\rho _{3},\cdots ,\rho
_{2n-1})\rightarrow H^{\ast }(PU(n))$
\end{quote}

\noindent \textsl{satisfies that}

\begin{enumerate}
\item[(6.5)] $\ker h=\langle \omega \otimes \theta (\xi _{I}),I\subseteq
Q_{p}^{+}(n)$\textsl{, }$p\in \left\{ p_{1},\cdots ,p_{t}\right\} \rangle $%
\textsl{.}
\end{enumerate}

\noindent \textbf{Proof.} Denote by $\varepsilon _{i}$ the cardinality of
the set $I_{i}$ and form the products

\begin{quote}
$b_{i}(I):=p_{1}^{\varepsilon _{1}}\cdots \widehat{p_{i}^{\varepsilon _{i}}}%
\cdots p_{t}^{\varepsilon _{t}}$, $1\leq i\leq t$,
\end{quote}

\noindent where $\widehat{}$ means omission. Since the sequence $%
\{b_{1}(I),\cdots ,b_{t}(I)\}$ so obtained is co-prime, there exists a
sequence $\{q_{1}(I),\cdots ,q_{t}(I)\}$ of integers such that

\begin{quote}
$q_{1}(I)\cdot b_{1}(I)+\cdots +q_{t}(I)\cdot b_{t}(I)=1$.
\end{quote}

\noindent Letting $J_{i}$ be the complement of $I_{i}\subseteq I$ and noting
that $b_{i}(I)\cdot \xi _{J_{i}}=c^{\ast }(\rho _{J_{i}})$ by (1.4), we get
by the additivity of $\theta $ the expansion

\begin{quote}
$\theta (\xi _{I})=\underset{1\leq i\leq t}{\Sigma }q_{i}(I)\cdot \theta
(b_{i}(I)\cdot \xi _{I})=\underset{1\leq i\leq t}{\Sigma }\pm q_{i}(I)\cdot
\theta (c^{\ast }(\rho _{J_{i}})\cup \xi _{I_{i}})$

$=\underset{1\leq i\leq t}{\Sigma }\pm q_{i}(I)\cdot \rho _{J_{i}}\cup
\theta (\xi _{I_{i}})$ (by the relation c) of Theorem B),
\end{quote}

\noindent where the sign $\pm $ depends on the reverse order of the
re-arrangement $J_{i}\sqcup I_{i}$ of the ordered sequence $I$. (6.4) is
shown by taking $a_{I}^{i}=\pm q_{i}(I)\cdot \rho _{J_{i}}$.\hfill $\square $

\bigskip

Assume that the integer $n$\ has the prime factorization $%
p_{1}^{r_{1}}\cdots p_{t}^{r_{t}}$. Our main result present the ring $%
H^{\ast }(PU(n))$ in the form of (6.1).

\bigskip

\noindent \textbf{Theorem C. }\textsl{The cohomology }$H^{\ast }(PU(n))$ 
\textsl{has the following decomposition into its subrings}

\begin{enumerate}
\item[(6.6)] $H^{\ast }(PU(n))=\Lambda (\rho _{3},\cdots ,\rho
_{2n-1})\oplus (\underset{1\leq i\leq t}{\oplus }\sigma _{p_{i}}(PU(n)))$%
\textsl{,}
\end{enumerate}

\noindent \textsl{where for each} \textsl{pair} $(p,r)=(p_{i},r_{i})$\textsl{%
\ with} $1\leq i\leq t$\textsl{,}

\begin{enumerate}
\item[(6.7)] $\sigma _{p}(PU(n))=\frac{J_{p,r}(\omega )\otimes \Lambda (\rho
_{3},\cdots ,\rho _{2n-1})}{\left\langle R_{I},\text{ }I\subseteq
Q_{p}^{+}(n)\QTR{sl}{,}l(I)\geq 2\right\rangle }$\textsl{\ with }$%
J_{p,r}(\omega )=\frac{\mathbb{Z}[\omega ]^{+}}{\left\langle p^{r}\omega
,p^{r-1}\omega ^{p},\cdots ,\omega ^{p^{r}}\right\rangle }$\textsl{,}
\end{enumerate}

\noindent \textsl{and where }

\begin{enumerate}
\item[(6.8)] $R_{I}=\underset{J\in S(I)}{\Sigma }p^{r-i_{1}-\varepsilon
_{I}(J)}\cdot \omega ^{\kappa _{I}(J)+1}\otimes \rho _{J}$ \textsl{(compare
with (6.3)).}
\end{enumerate}

\noindent \textbf{Proof. }Abbreviate the subring $\Lambda (\rho _{3},\cdots
,\rho _{2n-1})$ by $\Lambda _{n}$. By (5.6) the ring epimorphism $h:\mathbb{Z%
}[\omega ]\otimes \Lambda _{n}\rightarrow H^{\ast }(PU(n))$ (in the proof of
Theorem A) reduces to\textsl{\ }an epimorphism

\begin{quote}
$\widetilde{h}:J_{n}(\omega )\otimes \Lambda _{n}\rightarrow H^{\ast
}(PU(n)) $
\end{quote}

\noindent with $\ker \widetilde{h}=\langle \omega \cup \theta (\xi
_{I}),\left\vert I\right\vert \geq 2\rangle $. Moreover, if we put

\begin{quote}
$J_{p,r}(\omega ):=\frac{\mathbb{Z}[\omega ]^{+}}{\left\langle p^{r}\omega
,p^{r-1}\omega ^{p},p^{r-2}\omega ^{p^{2}},\cdots ,\omega
^{p^{r}}\right\rangle }$,
\end{quote}

\noindent then, by (5.7), the ring $J_{n}(\omega )\otimes \Lambda _{n}$ has
the decomposition into its free part and $p$-primary components:

\begin{enumerate}
\item[(6.9)] $J_{n}(\omega )\otimes \Lambda _{n}=\Lambda _{n}\oplus (%
\underset{1\leq i\leq t}{\oplus }J_{p_{i},r_{i}}(\omega )\otimes \Lambda
_{n})$.
\end{enumerate}

\noindent As a ring epimorphism, $\widetilde{h}$ must preserve the $p$%
-primary components. It follows that

\begin{quote}
$\sigma _{p}(PU(n))=0$ if $p\notin \{p_{1},\cdots ,p_{t}\}$,
\end{quote}

\noindent and that $\widetilde{h}$ restricts to an epimorphism

\begin{quote}
$\widetilde{h}_{i}:J_{p_{i},r_{i}}(\omega )\otimes \Lambda _{n}\rightarrow
\sigma _{p_{i}}(PU(n))$, $1\leq i\leq t$.
\end{quote}

On the other hand, by (6.4), each element $\omega \otimes \theta (\xi
_{I})\in \ker \widetilde{h}$ has the expansion

\begin{quote}
$a_{I}^{1}\cdot \omega \otimes \theta (\xi _{I_{1}})\oplus \cdots \oplus
a_{I}^{t}\cdot \omega \otimes \theta (\xi _{I_{t}})$ with $I_{i}\subset
Q_{p_{i}}^{+}(n)$,
\end{quote}

\noindent which is compatible with the decomposition (6.9) in the sense that

\begin{quote}
$a_{I}^{i}\cdot \omega \otimes \theta (\xi _{I_{i}})\in
J_{p_{i},r_{i}}(\omega )\otimes \Lambda _{n}$.
\end{quote}

\noindent This shows that $\widetilde{h}_{i}$ reduces furthermore to an
isomorphism

\begin{quote}
$\sigma _{p_{i}}(PU(n))\cong \frac{J_{p_{i},r_{i}}(\omega )\otimes \Lambda
_{n}}{\langle \omega \otimes \theta (\xi _{I}),I\subseteq
Q_{p_{i}}^{+}(n)\rangle }$, $1\leq i\leq t$.
\end{quote}

\noindent This is identical to (6.7), since by (6.3)

\begin{quote}
$\omega \otimes \theta (\xi _{I})=\frac{(-1)^{k-1}}{p^{r-i_{1}}}\binom{n}{%
p^{i_{1}}}\cdot R_{I}$,
\end{quote}

\noindent where the multiplier $\frac{(-1)^{k-1}}{p^{r-i_{1}}}\binom{n}{%
p^{i_{1}}}$ is co-prime to $p$ by Lemma 2.1.{\small \hfill $\square $}

\bigskip

\noindent \textbf{Remark 6.3. }Theorem C\textbf{\ }presents the cohomology $%
H^{\ast }(PU(n))$ as a ring. To clarify this we note that

i) Every summand in the right hand side of (6.6) is a subring of $H^{\ast
}(PU(n))$;

ii) The cup products between two different summands are subject to the
constraints:

\begin{quote}
$\cup :\Lambda (\rho _{3},\cdots ,\rho _{2n-1})\otimes \sigma
_{p_{i}}(PU(n))\rightarrow \sigma _{p_{i}}(PU(n))$;

$\cup :\sigma _{p}(PU(n))\otimes \sigma _{q}(PU(n))\rightarrow \{0\}$,
\end{quote}

\noindent where $p,q\in \{p_{1},\cdots ,p_{t}\}$ with $p\neq q$. In
addition, since $\Lambda (\rho _{3},\cdots ,\rho _{2n-1})$ has the additive
basis $\{1,\rho _{I}=\underset{j\in I}{\cup }{\small \ }\rho
_{2j-1},I\subseteq \left\{ 2,\cdots ,n\right\} \}$, the ideal $\sigma
_{p_{i}}(PU(n))$ is spanned additively by the elements

\begin{quote}
$\{\omega ^{r}\cdot \rho _{J},J\subseteq \left\{ 2,\cdots ,n\right\} ,1\leq
r<p^{r_{i}}\}$,
\end{quote}

\noindent and since $\rho _{2k-1}\in \Lambda (\rho _{3},\cdots ,\rho
_{2n-1}) $,

iii) The product between the summands $\Lambda (\rho _{3},\cdots ,\rho
_{2n-1})$ and $\sigma _{p_{i}}(PU(n))$ is given by

\begin{quote}
$\rho _{I}\cup \omega ^{r}\cdot \rho _{J}=0$ if $I\cap J\neq \emptyset $;

$\rho _{I}\cup \omega ^{r}\cdot \rho _{J}=\pm \omega ^{r}\cdot \rho
_{I\sqcup J}\func{mod}(\omega \otimes \theta (\xi _{I}),I\subseteq
Q_{p_{i}}^{+}(n))$ if $I\cap J=\emptyset $,
\end{quote}

\noindent where $I\sqcup J$ means the re-arrangement of the disjoint union $%
I\cup J$ into increasing order, and the sign $\pm $ depends on the reverse
order of the re-arrangement $I\sqcup J$ of $I\cup J$.{\small \hfill $\square 
$}

\bigskip

Theorem C is more convenient to use than Theorem A. We provide now such
evidences. Firstly, it is known that, for a compact and simply-connected Lie
group $G$ (e.g. \cite[Section \S 4.3]{D2}),

\begin{quote}
i) $\sigma _{p}(G)=0$ if $p\notin \{2,3,5\}$;

ii) $p\cdot \sigma _{p}(G)=0$ if $p\in \{2,3,5\}$.
\end{quote}

\noindent In contrast, Theorem C implies that

\bigskip

\noindent \textbf{Corollary 6.4.} \textsl{If the integer }$n$\textsl{\ has
the prime factorization }$p_{1}^{r_{1}}\cdots p_{t}^{r_{t}}$\textsl{, then}

\begin{quote}
\textsl{i) }$\sigma _{p}(PU(n))\neq 0$\textsl{\ if }$p\in \{p_{1},\cdots
,p_{t}\}$\textsl{;}

\textsl{ii) }$p^{s}\cdot \sigma _{p}(PU(n))\neq 0$\textsl{\ for all }$%
p=p_{i} $\textsl{\ and }$s<r_{i}$\textsl{.}\hfill $\square $
\end{quote}

\bigskip

Consider next the case $r_{1}=\cdots =r_{t}=1$. For a $p\in \{p_{1},\cdots
,p_{t}\}$ we have, in addition to

\begin{quote}
$J_{p,1}(\omega )=\frac{\mathbb{Z}_{p}[\omega ]^{+}}{\left\langle \omega
^{p}\right\rangle }$ and $Q_{p}^{+}(n)=\{1,p\}$,
\end{quote}

\noindent that $R_{\{1,p\}}=\omega \otimes \rho _{2p-1}$. Therefore, we get
from Theorem C that

\bigskip

\noindent \textbf{Corollary 6.5.} \textsl{If the integer }$n$\textsl{\ has
the prime factorization }$p_{1}\cdots p_{t}$\textsl{,} \textsl{then the ring 
}$H^{\ast }(PU(n))$ \textsl{is isomorphic to}

\begin{center}
$\Lambda (\rho _{3},\cdots ,\rho _{2n-1})\oplus (\underset{1\leq i\leq t}{%
\oplus }\frac{\mathbb{Z}_{p_{i}}[\omega ]^{+}}{\left\langle \omega
^{p_{i}}\right\rangle }\otimes \Lambda (\rho _{3},\cdots ,\widehat{\rho }%
_{2p_{i}-1},\cdots \rho _{2n-1}))$.\hfill $\square $
\end{center}

\bigskip

Suppose now that $n=p^{r}n^{\prime }$, where $p$ is a prime with $(n^{\prime
},p)=1$. Then $Q_{p}^{+}(n)=Q_{p}^{+}(p^{r})$. Since for each $p$-monotone $%
I\subseteq Q_{p}^{+}(n)$ the set $S(I)$, as well as the functions $%
\varepsilon _{I}$ and $\kappa _{I}$ on $S(I)$, relies only on the power $%
p^{r}$ (i.e. independent of $n^{\prime }$), we get by (6.6) the isomorphism

\begin{quote}
$\sigma _{p}(PU(n))=\sigma _{p}(PU(p^{r}))\otimes \Lambda (\rho
_{2p^{r}+1},\cdots ,\rho _{2n-1})$.
\end{quote}

\noindent This reduces the presentation of the ring $H^{\ast }(PU(n))$ to
that of $\sigma _{p}(PU(p^{r}))$.

\bigskip

\noindent \textbf{Corollary 6.6.} \textsl{If the integer }$n$\textsl{\ has
the prime factorization }$p_{1}^{r_{1}}\cdots p_{t}^{r_{t}}$\textsl{, then
the ring }$H^{\ast }(PU(n))$ \textsl{is isomorphic to}

\begin{center}
$\Lambda (\rho _{3},\cdots ,\rho _{2n-1})\oplus (\underset{1\leq i\leq t}{%
\oplus }\sigma _{p_{i}}(PU(p_{i}^{r_{i}}))\otimes \Lambda (\rho
_{2p_{i}^{r_{i}}+1},\cdots ,\rho _{2n-1}))$.\hfill $\square $
\end{center}

\bigskip

\noindent \textbf{Example 6.7.} If $n=2^{3}$ we have by Theorem C that

\begin{quote}
$H^{\ast }(PU(8))=\Lambda (\rho _{3},\cdots ,\rho _{17})\oplus \frac{%
J_{2,3}(\omega )\otimes \Lambda (\rho _{3},\cdots ,\rho _{17})}{\left\langle
R_{I},I\subseteq \{1,2,4,8\},l(I)\geq 2\right\rangle }$,
\end{quote}

\noindent where the generators $R_{I}$ of the ideal $\left\langle
R_{I},I\subseteq \{1,2,4,8\},l(I)\geq 2\right\rangle $ calculated by formula
(6.8) are tabulated below:

\begin{quote}
\begin{tabular}{l|l}
\hline\hline
$I\subseteq \{1,2,4,8\}$ & $R_{I}$ \\ \hline\hline
$\{1,2\}$ & $4\omega \otimes \rho _{3}$ \\ \hline
$\{1,4\}$ & $4\omega \otimes \rho _{7}+2\omega ^{2}\otimes \rho _{3}$ \\ 
\hline
$\{1,8\}$ & $4\omega \otimes \rho _{15}+2\omega ^{5}\otimes \rho _{7}+\omega
^{7}\otimes \rho _{3}$ \\ \hline
$\{2,4\}$ & $2\omega ^{2}\otimes \rho _{7}$ \\ \hline
$\{2,8\}$ & $2\omega ^{2}\otimes \rho _{15}+\omega ^{6}\otimes \rho _{7}$ \\ 
\hline
$\{4,8\}$ & $\omega ^{4}\otimes \rho _{15}$ \\ \hline
$\{1,2,4\}$ & $2\omega \otimes \rho _{3}\rho _{7}$ \\ \hline
$\{1,2,8\}$ & $2\omega \otimes \rho _{3}\rho _{15}+\omega ^{5}\otimes \rho
_{3}\rho _{7}$ \\ \hline
$\{1,4,8\}$ & $2\omega \otimes \rho _{7}\rho _{15}+\omega ^{3}\otimes \rho
_{3}\rho _{15}$ \\ \hline
$\{2,4,8\}$ & $\omega ^{2}\otimes \rho _{7}\rho _{15}$ \\ \hline
$\{1,2,4,8\}$ & $\omega \otimes \rho _{3}\rho _{7}\rho _{15}$ \\ \hline
\end{tabular}%
.\hfill \hfill $\square $
\end{quote}

\bigskip

\noindent \textbf{Remark 6.8. }For a compact and connected Lie group $G$
with a maximal torus $T$, the cohomology theory of the torus fibration (see
(3.1))

\begin{quote}
$T\rightarrow G\overset{\pi }{\rightarrow }G/T$
\end{quote}

\noindent connects two classical topics. The determination of the cohomology
ring $H^{\ast }(G/T)$ of the base manifold $G/T$ is a subject of the
Schubert calculus \cite[p.331]{Weil}, while the problem of computing the
cohomology ring $H^{\ast }(G)$ of the total space $G$ was raised by Cartan
in 1929 \cite{K}.

Based on the Schubert presentation of the ring $H^{\ast }(G/T)$ obtained in 
\cite{DZ}, we started in \cite{D2,DZ1} the project of constructing the
cohomology ring $H^{\ast }(G)$ in the context of Serre spectral sequence of $%
\pi $, where the project has been completed for the $1$-connected Lie groups 
$G$. The present work is devoted to extend the construction to the
non-simply connected groups $G=PU(n)$, in which the primary 1-forms $\rho
_{3}^{\prime },\cdots ,\rho _{2n-1}^{\prime }$ constructed in Definition 4.2
play a central role.\hfill \hfill $\square $

\bigskip

Haibao Duan

dhb@math.ac.cn

Yau Mathematical Science Center, Tsinghua University, Beijing 100084;

Academy of Mathematics and Systems Sciences, Chinese Academy of Sciences,
Beijing 100190.

\end{document}